\documentclass[a4paper, 11pt]{article} 
\usepackage[T1]{fontenc}

\usepackage[margin={1.25cm,1.75cm},footskip=.75cm]{geometry}
\usepackage[margin=1.25cm,font={footnotesize}]{caption}

\usepackage[colorlinks,citecolor=blue,linkcolor=blue,urlcolor=blue]{hyperref}
\usepackage{amsmath}
\usepackage{amssymb}
\usepackage{bm}
\usepackage{graphicx}
\usepackage{booktabs}
\usepackage{appendix}
\usepackage{doi}
\usepackage{hyphenat}
\usepackage{float}
\usepackage[authoryear,round]{natbib}
\usepackage{stmaryrd} 
\usepackage{mathtools} 
\usepackage{subfigure}
\usepackage[affil-it]{authblk}
\usepackage{fancyhdr}
\hypersetup{pdfcreator={LNoel},pdfproducer={}}

\title{\bf XIGA: An eXtended IsoGeometric Analysis approach for multi-material problems}
\author{L. No{\"e}l$^{1}$, M. Schmidt$^{2}$, K. Doble$^{2}$, J.A. Evans$^{2}$, and K. Maute$^{2}$}
\affil{\small
$^1${Department of Precision and Microsystems Engineering, Faculty of Mechanical, Maritime and Materials Engineering,}\\ 
{Delft University of Technology}, {The Netherlands}\\
\href{mail to:l.f.p.noel@tudelft.nl}{l.f.p.noel@tudelft.nl}\\
$^2${Aerospace Mechanics Research Center, Department of Aerospace Engineering Sciences},\\
{University of Colorado Boulder}, {CO}, {USA}\\
\href{mail to:k.maute@colorado.edu}{k.maute@colorado.edu}\\
}
\date{}
\providecommand{\keywords}[1]{\textbf{{Keywords --}} #1}

\begin{document}

\maketitle

\abstract{Multi-material problems often exhibit complex geometries along with physical responses presenting large spatial gradients or discontinuities. In these cases, providing high-quality body-fitted  finite element analysis meshes and obtaining accurate solutions remain challenging. Immersed boundary techniques provide elegant solutions for such problems. Enrichment methods alleviate the need for generating conforming analysis grids by capturing discontinuities within mesh elements. Additionally, increased accuracy of physical responses and geometry description can be achieved with higher-order approximation bases. In particular, using B-splines has become popular with the development of IsoGeometric Analysis. In this work, an eXtended IsoGeometric Analysis (XIGA) approach is proposed for multi-material problems. The computational domain geometry is described implicitly by level set functions. A novel generalized Heaviside enrichment st\-rategy is employed to accommodate an arbitrary number of materials without artificially stiffening the physical response. Higher-order B-spline functions are used for both geometry representation and analysis. Boundary and interface conditions are enforced weakly via Nitsche's meth\-od, and a new face-oriented ghost stabilization meth\-odology is used to mitigate numerical instabilities arising from small material integration subdomains. Two- and three-dimensional heat transfer and elasticity problems are solved to validate the approach. Numerical studies provide insight into the ability to handle multiple materials considering sharp-edged and curved interfaces, as well as the impact of higher-order bases and stabilization on the solution accuracy and conditioning.}\\


\section{Introduction}\label{sectionIntroduction}
Multi-material problems play an important role for a wide range of applications in engineering, such as problems involving composite or functionally graded materials, multiple phase interactions, or contact between components. However, efficiently and accurately predicting the physical response described by partial differential equations of multi-material problems remains challenging. Such problems often exhibit complex geometries, intricate material arrangements, or small features along with physical responses presenting large spatial gradients or discontinuities. Therefore, providing an accurate resolution for both the geometry and the physics around material interfaces is crucial. When performing analysis with traditional finite element methods (FEM), this can only be achieved by generating highly-refined body-fitted approximation meshes, which is a tedious task and is known to represent a substantial part of the analysis time, see \cite{2010BazilevsEtAl}. This is especially the case when the number of materials increases and when boundaries and material interfaces change in time due to physical, e.g., time evolving interfaces, see \cite{2015KamenskyEtAl}, or numerical processes, e.g., topology optimization, see \cite{2020NoelEtAl}.

Over the past decades, immersed boundary techniques have gained in popularity, providing an elegant solution for the multi-material problems characterized above. These methods can accommodate complex domain boundaries and interfaces without the need to construct a conforming body-fitted mesh. The first immersed boundary method was formally introduced by \cite{1972Peskin}. Since this first occurrence, several immersed finite element approaches have been presented in the literature. Fictitious domain methods, also known as embedded domain methods, circumvent the need to generate conforming analysis meshes by embedding the computational domain in a larger one and applying specific integration techniques. Variants of the method are based on different approaches to impose boundary and interface conditions with penalty methods by \cite{2007RamiereEtAl}, with Lagrange multipliers by \cite{1994GlowinskiEtAl}, \cite{2007GlowinskiKuznetsov}, or \cite{2010BurmanHansbo}, or with Nitsche's method by \cite{2002HansboHansbo}, \cite{2009DolbowHarari}, \cite{2012BurmanHansbo}, and \cite{2015BurmanEtAl}. Conversely, enrichment methods alleviate the need for generating conforming approximation meshes by capturing \textit{a priori} known discontinuous behaviors within the mesh elements. Originally developed to represent moving fronts and crack propagation, enrichment based approaches have been extended to tackle various types of interface problems with strong and weak discontinuities. Among these techniques, a few noticeable ones are the Partition of Unity Method (PUM) proposed by \cite{1997BabuskaMelenk}, the Generalized Extended Finite Element Method (GFEM) as introduced in \cite{2000aStrouboulisEtAl,2000bStrouboulisEtAl}, the eXtended Finite Element Method (XFEM) as proposed by \cite{1999MoesEtAl} and \cite{1999BelytschkoBlack}, and the Interface enriched Generalized Finite Element Method (IGFEM) introduced by \cite{2012SoghratiEtAl}.

In most of the aforementioned papers, the implementation of immersed boundary techniques relies on low order approximation functions, in particular linear Lagrange basis functions, to represent both the geometry and the physics. Such a choice of basis functions suffers from several shortcomings in terms of geometry resolution, and accuracy of physical responses. Using \textit{p}-version FEM, see \cite{1988Babuska} for FEM based on higher-order Lagrange functions or \cite{1984Patera} and \cite{1999KarniadakisSherwin} for spectral FEM based on higher-order spectral basis functions, enables improved accuracy of physical responses per degree of freedom (DOF) and leads to higher convergence rates with mesh refinement, i.e., \textit{h}-refinement. Additionally, using higher-order functions for the geometry representation improves accuracy in the presence of curved interfaces and boundaries.

Several papers in the literature use higher-order ba\-ses, such as Lagrange or spectral ones, in combination with immersed boundary techniques. \cite{2007ParvizianEtAl} and \cite{2008DusterEtAl} introduced the finite cell meth\-od. Similar to other fictitious domain approaches, the method extends the analysis domain to embed the physical one, but makes use of higher-order \textit{Ansatz} functions to approximate the extended variables. Numerous contributions focused on enrichment methods to accurately represent the geometry of and the physics around curved cracks. \cite{2002WellsEtAl} studied the propagation of displacement discontinuities in strain-softening media with second order Lagrange bases. Working on crack propagation, \cite{2003StaziEtAl} used second order Lagrange bases for the finite element approximations, while \cite{2003ZiBelytschko} extended this approach to higher-order enrichment functions. Tackling both strong and weak discontinuities, \cite{2009ChengFries} resolved curved boundaries and interfaces by generating integration subcells with one curved side and by applying corrections to the enrichment formulation. Focusing on material interfaces, \cite{2010Dreau} exploited the XFEM with a corrected enrichment and represented the geometry on sub-meshes finer than the one used for the mechanical fields. A similar approach was proposed in \cite{2012LegrainEtAl}. \cite{2011HaasemannEtAl} proposed a numerical integration strategy based on NURBS surfaces for higher-order XFEM and weak discontinuities. \cite{2016Lehrenfeld} used parametric mappings of the integration cells to reduce the interface representation error and the associated integration error. More recently, \cite{2020SaxbyHazel} proposed a higher-order modified XFEM based on corrected basis functions for weak discontinuity problems. Wor\-king with higher-order spectral basis functions, \cite{2005LegayEtAl} and more recently \cite{2019ChinSukumar} proposed a spectral XFEM approach to tackle weak discontinuity problems with curved interfaces.
	
Along with the development of IsoGeometric Analysis (IGA), using B-splines or NURBS as basis functions has become an increasingly popular approach, see \cite{2005HuguesEtAl} and \cite{2009CottrellEtAl}. In IGA, both the geometry of a structure and its physical behavior are described using splines. From a geometry point of view, using B-splines and NURBS facilitates compatibility with Computer Aided Design (CAD) software. From an analysis point of view, using smooth and higher-order bases, such as quadratic and cubic B-splines, leads to more accurate physical responses per DOF than traditional $C^0$ finite element approaches, see \cite{2008HuguesEtAl}; \cite{2009EvansEtAl}; \cite{2014HuguesEtAl}.

Over the years, several research works have aimed at combining the advantages of immersed boundary techniques and smooth higher-order basis functions, such as B-splines, NURBS, and other variants. Based on fictitious domain approaches, \cite{2012SchillingerEtAl} proposed a B-spline version of the finite cell method. \cite{2015KamenskyEtAl} further extended the concept to tackle fluid-solid interaction problems and introduced the term \textit{immersogeometric analysis}. \cite{2001HolligEtAl} introduced the web-method using weighted extended B-splines as basis functions to solve Dirichlet problems. Modeling weak discontinuities and in particular material interfaces, \cite{2011SanchesEtAl} developed an immersed boundary technique based on B-spline bases and a modified basis to locally interpolate the Dirichlet boundary conditions. Focusing on enrichment methods, \cite{2015JiangEtAl} proposed a robust Nitsche's method to tackle interface problems with the XFEM based on B-spline basis functions. They used a separate locally refined mesh to improve the geometry representation. \cite{2015JiaEtAl} solved curved material interface problems with XFEM based on NURBS and used curved integration elements for increased accuracy. To further resolve the interface geometry, \cite{2019ChenEtAl} implemented the XFEM with locally refined B-splines to allow for adaptive local refinement around the interfaces. Recently, \cite{2018ElfversonEtAl,2019ElfversonEtAl} proposed a so-called cutIGA approach and a symmetric Nitsche's method for imposing boundary conditions, as well as a drop of basis functions for improved stability.

To date, the scope and the applications of immersed boundary techniques with higher-order basis functions are rather limited in terms of the number of phases or materials, but also in terms of the geometric complexity considered. Tools to handle the geometric representation were investigated in numerous publications. \cite{2011TranEtAl} used several level set functions (LSF) to accurately represent complex microstructures with multiple spatially close inclusions and avoid numerical artefacts using the XFEM. \cite{2011MoumnassiEtAl} also used several LSFs to accurately represent sharp features and curved interface without mesh refinement. \cite{2011XiaEtAl} proposed a matched interface and boundary method to tackle multi-material and triple junctions. \cite{2012HouEtAl} built specific approximations for interface elements presenting multiple material and triple-junctions. A similar approach has recently been proposed by \cite{2020ChenEtAl}, who treated triple-junction points with an immersed boundary technique through the construction of specific functions on interface elements. \cite{2014Soghrati} extended the IGFEM to handle multi-material interfaces by constructing special enrichment functions. However, most proposed frameworks lack versatility and do not offer a systematic approach to tackle multi-material problems in two and three dimensions.

In this paper, we propose a versatile XIGA approach to tackle multi-material problems in two and three dimensions. The geometry of the computational domain is represented implicitly by one or multiple LSFs. The LSFs are used to determine subregions of the computational domain that are associated to different phases and materials. This approach allows for handling of straight-edged and curved interfaces, as well as \textit{N}-ma\-terial junctions, in a systematic way. The governing equations are integrated separately on each material subdomain, and elements where multiple materials coexist are decomposed into single material integration subdomains. The finite element approximations for both the geometry and the mechanical fields use multi-variate B-splines. They are smooth higher-order basis functions and provide higher accuracy per DOF and higher convergence rates than traditional $C^0$ finite element bases. In this paper, a novel generalized Heaviside enrichment strategy is used with multiple enrichment levels to introduce discontinuities at external boundaries and material interfaces. Boundary and interface conditions are weakly enforced using Nitsche's method. Numerical instabilities associated with small material integration subdomains are mitigated by an adapted version of the face-oriented ghost stabilization. The combination of these ingredients results in a versatile and robust approach to tackle multi-material problems.

In most enrichment approaches, the approximation space is extended using different enrichments for each material domain. However, this approach can lead to an artificially stiffened physical response when one or more material domains are disconnected. To alleviate this issue, the approximation space was extended using different enrichments for each connected material subdomain in \cite{2003TeradaEtAl}, or \cite{2004HansboHansbo}. Nonetheless, this approach can still lead to an artificially stiffened physical response when the intersection of a connected material subdomain with the support of a particular basis background function is disconnected, see \cite{2014MakhijaMaute}. This situation frequently arises when B-spline basis functions are employed rather than classical finite element basis functions. This is because B-spline basis functions have larger support regions than classical finite element basis functions. To tackle this issue, we build on previous work on Lagrange basis functions (see \cite{2014MakhijaMaute}) and enrich each individual basis function separately based on the topology of the material layout within the basis function support.

For multi-material problems, ghost stabilization procedures typically rely on the polynomial extension of function values within a material region of an element adjacent to a ghost facet to the ghost facet itself, see \cite{2014BurmanHansbo}. However, with the enrichment strategy proposed in this work, such an extension is not well-defined when a material region within the element is disconnected. In particular, the polynomial extensions associated with different connected material subregions may differ. To overcome this issue, we introduce a ghost stabilization strategy that explicitly accounts for the topology of the material layout of elements adjacent to ghost facets.

The remainder of the paper is organized as follows. Section \ref{sectionGeometry} focuses on the use of one or multiple LSFs to represent the geometry of external boundaries and material interfaces. The proposed XIGA formulation is detailed in Section \ref{sectionXIGA}. First, a brief discussion of B-splines as basis functions for finite element analysis is provided. Then, the immersed boundary technique, i.e., here the XFEM, is detailed in terms of enrichment strategy, creation of the model, formulation of the governing equations, stabilization via an adapted face-oriented ghost stabilization, and integration. Section \ref{sectionExamples} illustrates the capabilities and the robustness of the proposed XIGA approach by solving canonical two- and three-dimensional problems focusing on heat conduction and elasticity. Finally, Section \ref{sectionConclusion} draws conclusions about the developed XIGA approach and proposes directions for future work.

\section{Geometry representation}\label{sectionGeometry}
Although the proposed XIGA approach is not restricted to any particular geometry description method, the geometry of a computational domain is represented by one or multiple LSFs in this paper. This specific geometry representation was chosen to ease future work tackling evolving interfaces, in particular for design through topology optimization. 

The level set method (LSM) was developed by \cite{1988OsherSethian} to efficiently track front propagation. The method allows for the implicit representation of a geometry by a LSF. An iso-level $\phi_t$ of the LSF, generally chosen equal to 0, describes the interface $\Gamma_{\pm}$ between two regions $\Omega_+$ and $\Omega_-$ of the analysis domain $\Omega$ via:
\begin{equation}
\begin{array}{ll}
\phi(\mathbf{x}) < \phi_t, & \forall\, \mathbf{x} \in \Omega_+,\\
\phi(\mathbf{x}) > \phi_t, & \forall\, \mathbf{x} \in \Omega_-,\\
\phi(\mathbf{x}) = \phi_t, & \forall\, \mathbf{x} \in \Gamma_{\pm}.\\
\end{array}
\label{eqLevelSet}
\end{equation}
\begin{figure*}[!h]\center
\includegraphics[scale=1]{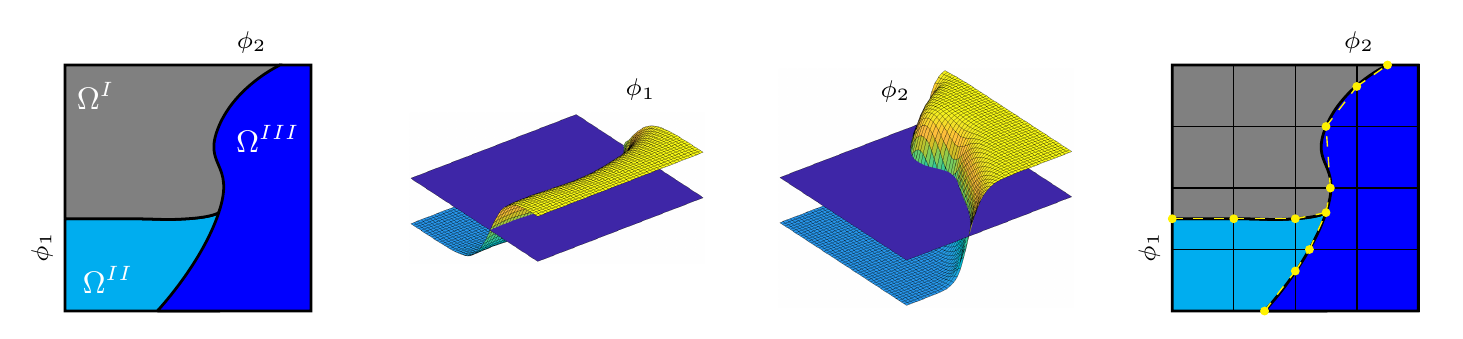}
\caption{Geometry description of a computational domain, made of three material domains $\Omega^{I}$, $\Omega^{II}$, and $\Omega^{III}$, using two LSFs, $\phi_1$ and $\phi_2$.}
\label{figGeometryDescription}
\end{figure*}

An example of the type of multi-material problems addressed in this work is given in Fig. \ref{figGeometryDescription}, where a computational domain $\Omega$, made of three material domains $\Omega^{I}$, $\Omega^{II}$, and $\Omega^{III}$, is described using two LSFs, $\phi_1$ and $\phi_2$. In this paper, multiple LSFs $\phi_i(\mathbf{x})$, $i=1, \dots, n$, are used to describe the external boundaries of and the material interfaces within the computational domain $\Omega$. A multi-phase level set model is exploited, as introduced by \cite{2002VeseChan}. With $n$ LSFs, a maximum of $2^n$ phases can be represented. In this paper, a phase represents a subregion of the analysis domain associated with a unique combination of positive or negative valued LSFs. A phase index $\mathcal{P}$ is assigned to each subregion based on the LSF signs. The phase assignment procedure is sequenced as follows. First, characteristic functions $f_i(\mathbf{x})$, $i=1, \dots, n$, are used to characterize the point $\mathbf{x}$ with respect to an iso-level $\phi_t$ of the LSF $\phi_i(\mathbf{x})$, i.e., whether $\mathbf{x}$ is inside, outside, or on the iso-level contour, as:  
\begin{equation}
f_i(\mathbf{x}) = \left\{
\begin{array}{ll}
0, & \phi_i(\mathbf{x}) < \phi_t\\
1, & \phi_i(\mathbf{x}) > \phi_t\\
\mbox{on interface,} & \phi_i(\mathbf{x}) = \phi_t\\
\end{array}
\right.
\label{eqInsideOutside}
\end{equation}

These characteristic functions $f_i(\mathbf{x})$ are used to assign a unique index $\mathcal{P}(\mathbf{x})$ to the point $\mathbf{x}$:
\begin{equation}
\mathcal{P}(\mathbf{x}) = \sum_{j=1}^{n} 2^{j-1}\ f_{j}(\mathbf{x}).
\label{eqPhaseIndex}
\end{equation}

Finally, a material describing the constitutive behavior is assigned to each phase. The phase indices $\mathcal{P}(\mathbf{x})$ are associated with the corresponding material indices $\mathcal{M}(\mathbf{x})$ through a map $m$ following:
\begin{equation}
\mathcal{M}(\mathbf{x}) = m(\mathcal{P}(\mathbf{x})).
\label{eqPhaseAssignment}
\end{equation}

\begin{figure*}[!h]\center
\includegraphics[scale=1]{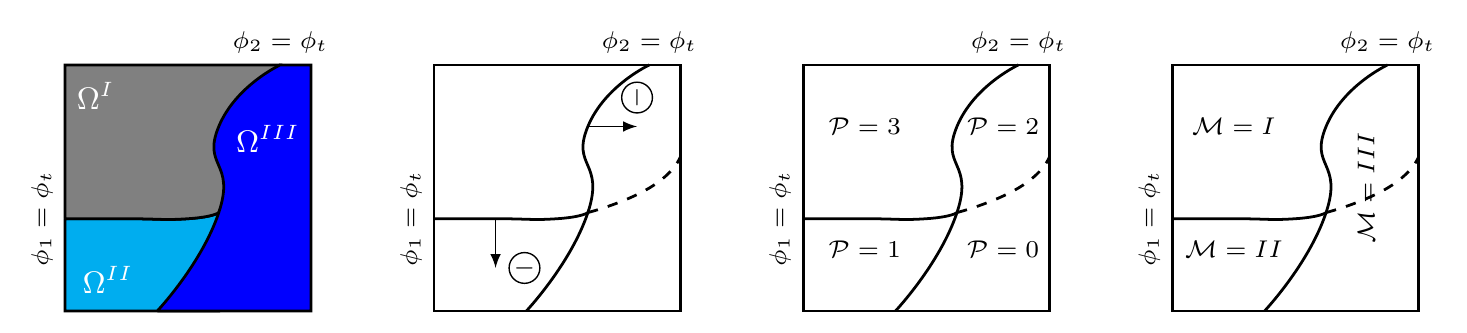}
\caption{Phase and material assignment procedure for a three-material problem described by two LSFs, $\phi_1$ and $\phi_2$.}
\label{figPhaseAssignment}
\end{figure*}

The phase and material assignment procedure is illustrated with the three-material problem in Fig. \ref{figPhaseAssignment}. First, the characteristic functions $f_i(\mathbf{x})$ are evaluated based on the LSFs signs. In this picture, the minus sign indicates that $\phi_i(\mathbf{x}) < \phi_t$ and thus, $f_i(\mathbf{x}) = 0$. Then, the phase indices $\mathcal{P}(\mathbf{x})$ are computed by Eq.\eqref{eqPhaseIndex} based on the characteristic functions $f_i(\mathbf{x})$. Finally, a material is assigned to each phase, by associating a material index $\mathcal{M}$ to each phase index. 

Each LSF $\phi_i(\mathbf{x})$ is discretized on a mesh using basis functions $B_k(\mathbf{x})$ as:
\begin{equation}
\phi_i^{h}(\mathbf{x}) = \sum_{k} B_k(\mathbf{x})\, \phi_{i}^{k},
\label{eqLevelSetDiscretization}
\end{equation}
where $\phi_{i}^{k}$ are the coefficients associated to the LSF $\phi_i(\mathbf{x})$. In this work, B-spline basis functions are chosen. The coefficients and corresponding basis functions are used to evaluate nodal level set values on the background mesh generated for analysis, see Section \ref{sectionXIGA}. The LSF is interpolated linearly along the element edges to determine the intersection of the $\phi=\phi_t$ iso-contour with the element edges, see Subsection \ref{subsectionIntegration}. Using a linear interpolation inherently leads to a low order approximation of geometry, which might limit the accuracy of the physical response analysis. This issue can be alleviated by first refining the background mesh and then interpolating the LSF on this refined mesh.

\section{XIGA formulation}\label{sectionXIGA}
This section focuses on the XIGA approach proposed in this paper. First, B-spline basis functions are briefly reviewed in Subsection \ref{subsectionB-spline}. Then, our novel enrichment strategy based on generalized Heaviside functions, used to accommodate multiple materials within a basis function support, is explained in Subsection \ref{subsectionEnrichment}. Subsection \ref{subsectionGoverning} summarizes the formulation of the governing equations considered in this work for heat conduction and elasticity. The techniques to enforce boundary and interface conditions and to stabilize the XIGA formulation are detailed in Subsections \ref{subsectionBoundary} and \ref{subsectionStabilization} respectively. As several materials may coexist within the same element, special attention is required to perform numerical integration. This integration procedure is described in Subsection \ref{subsectionIntegration}.

\subsection{B-splines for finite element analysis}\label{subsectionB-spline}
In this paper, B-splines are used to approximate the level set and physics variable fields. This subsection briefly recalls basic concepts of constructing B-splines in one and multiple dimensions. 

Considering a knot vector $\Xi = \{ \xi_{1},\xi_{2}, \dots ,\xi_{n+p+1} \}$, for which $\xi \in \mathbb{R}$ and $\xi_{1} \leq \xi_{2} \leq \dots \leq \xi_{n+p+1}$, a univariate B-spline basis function $N_{i,p}(\xi)$ of degree $p$ is constructed recursively starting from the piecewise constant basis function:
\begin{equation}
N_{i,0}(\xi) =
\begin{cases}
1, & \text{if}\ \xi_{i} \leq \xi \leq\xi_{i+1},\\
0, & \text{otherwise}.
\end{cases}
\end{equation}
The Cox de Boor recursion formula is used to obtain the basis functions for higher degrees $p>0$, see \cite{1972DeBoor}:
\begin{equation}
N_{i,p}(\xi) = \frac{\xi - \xi_{i}}{\xi_{i+p} - \xi{i}}\ N_{i,p-1}(\xi)
+ \frac{\xi_{i+p+1} - \xi}{\xi_{i+p+1} - \xi_{i+1}}\ N_{i+1,p-1}(\xi).
\label{eqUnivariateB-spline}
\end{equation}
A knot is said to have a multiplicity $k$ if it is repeated $k$ times in the knot vector. The corresponding B-spline basis exhibits a $C^{p-k}$ continuity at that specific knot, while it is $C^{\infty}$ in between unique knots.

To tackle \textit{n}-dimensional problems, multi-variate B-spline basis functions $B_{\mathbf{i},\mathbf{p}}(\boldsymbol{\xi})$ are obtained by the tensor product of univariate B-spline basis functions. Denoting the parametric space dimension by $d_{p}$, a tensor-product B-spline basis is constructed starting from $d_{p}$ knot vectors $\Xi^{m} = \{\xi_{1}^{m},\xi_{2}^{m},\-\dots,\-\xi_{n_{m}+p_{m}+1}^{m}\}$ with $p_{m}$ being the polynomial degree and $n_{m}$ the number of basis functions in the parametric direction $m = 1, \dots, d_{p}$. A tensor-product B-spline basis function is generated from $d_{p}$ univariate B-splines $N_{i_{m},p_{m}}^{m}(\xi^{m})$ in each parametric direction $m$ using the formula:
\begin{equation}
B_{\mathbf{i},\mathbf{p}}(\boldsymbol{\xi}) = \prod_{m=1}^{d_{p}} N_{i_{m},p_{m}}^{m}(\xi^{m}),
\label{Eq_multiTensorProd}
\end{equation}
where the position in the tensor product structure is given by the index $\mathbf{i} =\{ i_{1}, \dots, i_{d_{p}} \}$, and the polynomial degree is denoted by $\mathbf{p} = \{ p_{1}, \dots, p_{d_{p}} \}$.

In this paper, Lagrange extraction, as introduced by \cite{2016SchillingerEtAl}, is used to facilitate a classical finite element implementation of the integration of the governing equations over the mesh of background elements, here defined as the tensor product of nonempty knot spans. This approach avoids the need to consider the non-elemental-locality of B-splines and simplifies the integration procedure for elements occupied by multiple materials.

\subsection{Enrichment strategy}\label{subsectionEnrichment}
The XFEM was introduced by \cite{1999MoesEtAl} and \cite{1999BelytschkoBlack} to model crack propagation without remeshing. The method enables the prediction of discontinuous or singular behaviors within an element by adding specific enrichment functions to the classical finite element approximation. In this paper, we follow the work by \cite{2003TeradaEtAl}, \cite{2004HansboHansbo}, and \cite{2014MakhijaMaute}, and use a generalized Heaviside enrichment strategy to introduce discontinuities along geometries and material interfaces. We further generalize the approach and enrich each basis function separately based on the material layout within the basis function support to ensure independent approximation on each connected material subregion. 

\begin{figure*}[!h]\center
\includegraphics[scale=1]{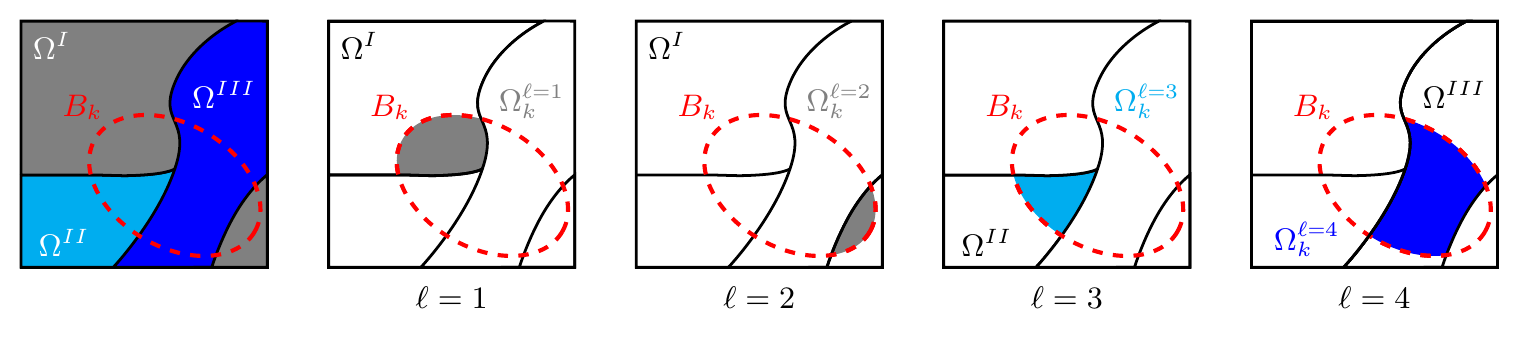}
\caption{Enrichment strategy for a basis function $B_k$ spanning the three material subdomains in a three-material problem.}
\label{figEnrichmentStrategy}
\end{figure*}

The enrichment level selection procedure proposed in this paper is illustrated for a three-material problem in Fig. \ref{figEnrichmentStrategy} for a basis function $B_k$ spanning the three material subdomains. The basis function support is delimited by a red dashed line. Within this support, $\mbox{supp}(B_k)$, four separate connected material subregions $\Omega_k^{\ell}$ exist, each occupied by one and only one material, such that $\mbox{supp}(B_k) = \bigcup_{\ell=1}^{4} \Omega_k^{\ell}$. Two subregions $\Omega_k^{\ell=1}$ and $\Omega_k^{\ell=2}$ are occupied by the same material $I$. Only one subregion $\Omega_k^{\ell=3}$ is occupied by material $II$ and one $\Omega_k^{\ell=4}$ is occupied by material $III$. As four material subregions $\Omega_k^{\ell}$ exist within the basis support, four enrichment levels $L_k=4$ are necessary.

Formally, considering a multi-material problem, a state field $\mathbf{u}(\mathbf{x})$ is approximated as: 
\begin{equation}
\mathbf{u}^h(\mathbf{x}) = \sum_{k=1}^{K} \sum_{\ell=1}^{L_k} \varphi_k^{\ell}(\mathbf{x}) B_k(\mathbf{x})\ u_{k}^{\ell},
\end{equation}
where $K$ is the number of background basis functions, and $L_k$ is the number of separate connected material subregions $\{ \Omega_k^{\ell} \}_{\ell=1}^{L_k}$ in the support of background basis function $B_k$ for $k=1,\dots,K$. The coefficient $u^{\ell}_{k}$ is the DOF associated with background basis function $B_k$ and material subregion $\Omega_{k}^{\ell}$ for $k=1,\dots,K$ and $\ell=1,\dots,L_k$. The function $\varphi_k^{\ell}(\mathbf{x})$ is an indicator function that determines whether a point $\mathbf{x}$ belongs to a material subregion $\Omega_k^{\ell}$:
\begin{equation}
\varphi_k^{\ell}(\mathbf{x}) = I_{\Omega_k^{\ell}}(\mathbf{x}) = \left\{
\begin{array}{ll}
1 & \displaystyle \mbox{if}\ \mathbf{x} \in \Omega_k^{\ell},\\[5pt]
0 & \displaystyle \mbox{otherwise}.
\end{array}
\right.
\label{Eq_enrichment_selectionFunction}
\end{equation}
The set of functions $\mathcal{B}:= \left\{ \varphi_k^{\ell} B_k : k \in \{ 1,\dots,K\}\ \mbox{and}\ \ell \in \{ 1,\dots,L_k \}  \right\}$ possesses several useful properties. First, note that for each $\mathbf{x}\in \mbox{supp}(N_k)$, there is one and only one $\ell$ for which $\varphi_k^{\ell}(\mathbf{x})\neq 0$. Thus, $\sum_{\ell} \varphi_k^{\ell}(\mathbf{x}) B_k(\mathbf{x}) = B_k(\mathbf{x})$, and as the background basis $\{ B_k \}_{k=1}^{K}$ forms a partition of unity, it follows that the functions in $\mathcal{B}$ do as well. The functions in $\mathcal{B}$ are also pointwise non-negative. Finally, as the background basis functions are locally linearly independent, so are the functions in $\mathcal{B}$. This indicates the functions in $\mathcal{B}$ form a basis. We refer to this basis as the enriched basis as it derives from the enrichment strategy.

\subsection{Governing equations in discretized form}\label{subsectionGoverning}
The proposed XIGA approach is not limited to any particular type of partial differential equation. However, in this paper, we restrict our attention to elliptic problems, namely multi-material linear elasticity and heat conduction problems. Either a linear elastic or linear diffusive material is assumed for each non-void domain.

In this work, the total residual $\mathcal{R}$, i.e., the discrete form of the governing equations, consists of four terms which are discussed subsequently. We solve for static equilibrium to enforce balance of linear momentum with\-in each material domain $\Omega^I$, where $I$ is the material index:
\begin{equation}
\mathcal{R}(\mathbf{u},\delta \mathbf{u}) = \mathcal{R}_{Lin}^{\mathbf{u}} + \mathcal{R}_{D}^{\mathbf{u}} + \mathcal{R}_{Itf}^{\mathbf{u}} + \mathcal{R}_{Ghost}^{\mathbf{u}} = 0,
\label{Eq_RFullU}
\end{equation}
where $\mathbf{u}$ and $\delta \mathbf{u}$ are the displacement field and the test function, respectively.

We solve for static equilibrium to enforce heat balance within each material domain $\Omega^I$, where $I$ is the material index:
\begin{equation}
\mathcal{R}(\theta,\delta \theta) = \mathcal{R}_{Lin}^{\theta} + \mathcal{R}_{D}^{\theta} + \mathcal{R}_{Itf}^{\theta} + \mathcal{R}_{Ghost}^{\theta}  = 0,
\label{Eq_RFullTheta}
\end{equation}
where $\theta$ and $\delta \theta$ are the temperature field and the test function, respectively.

The first residual term $\mathcal{R}^{\mathbf{u}}_{Lin}$ for linear elasticity reads:
\begin{equation}
\displaystyle
\mathcal{R}^{\mathbf{u}}_{Lin} = \sum_I \Bigg[\
\textcolor{white}{+} \int_{\Omega^I}\ \delta \boldsymbol{\varepsilon} : \boldsymbol{\sigma}\ d\Omega 
+ \int_{\Omega^I}\ \delta \mathbf{u} \cdot \mathbf{b}\ d\Omega 
- \int_{\Gamma_N^I} \delta \mathbf{u} \cdot \mathbf{t}_{N}\ d\Gamma\ \Bigg],
\label{Eq_RLinU}
\end{equation}
where body loads, $\mathbf{b}$, are acting on the domain $\Omega^I$ and traction forces, $\mathbf{t}_N$, are applied on the Neumann boundary, $\Gamma_N^I$. The Cauchy stress tensor is denoted by $\boldsymbol{\sigma} = \mathbf{D}\ \boldsymbol{\varepsilon}$ and is obtained by multiplication of the infinitesimal strain tensor $\boldsymbol{\varepsilon} = \frac{1}{2} \left( \nabla \mathbf{u} + \nabla \mathbf{u}^T \right)$ with the fourth order constitutive tensor $\mathbf{D}$, here for isotropic linear elasticity, expressed as a function of the Young's modulus $E$ and the Poisson ratio $\nu$ in Voigt notation in 2D as:
 \begin{equation} 
\mathbf{D} = \tilde{E}\left[
\begin{array}{cccccc}
1-\nu & \nu & 0\\
\nu & 1-\nu & 0\\
0 & 0 & \frac{1-2\nu}{2}\\
\end{array}
\right],
\label{Eq_Dmatrix_2D}
\end{equation} 
and in 3D as:
\begin{equation}
\mathbf{D} = \tilde{E}\left[
\begin{array}{cccccc}
1-\nu & \nu & \nu & 0 & 0 & 0\\
\nu & 1-\nu & \nu & 0 & 0 & 0\\
\nu & \nu & 1-\nu & 0 & 0 & 0\\
0 & 0 & 0 & \frac{1-2\nu}{2} & 0 & 0\\
0 & 0 & 0 & 0 & \frac{1-2\nu}{2} & 0\\
0 & 0 & 0 & 0 & 0 & \frac{1-2\nu}{2}\\
\end{array}
\right],
\label{Eq_Dmatrix_3D}
\end{equation}
with $\displaystyle \tilde{E} = \frac{E}{(1+\nu)(1-2\nu)}$.

The first residual term $\mathcal{R}^{\theta}_{Lin}$ for heat conduction reads:
\begin{equation}
\displaystyle
\mathcal{R}^{\theta}_{Lin} = \sum_I \Bigg[\
- \int_{\Omega^I}\ \delta \nabla \theta \cdot \mathbf{q}\ d\Omega
+ \int_{\Omega^I}\ \delta \theta\ q_B\ d\Omega
- \int_{\Gamma_N^I} \delta \theta\ q_{N}\ d\Gamma\ \Bigg],
\label{Eq_RLinTheta}
\end{equation}
where body heat loads, $q_B$, are acting on the domain $\Omega^I$ and heat fluxes, $q_N$, are applied on the Neumann boundary, $\Gamma_N^I$. The heat flux, $\mathbf{q} = - \left( \boldsymbol{\kappa} \cdot \nabla \theta \right)$, is obtained by multiplying the conductivity tensor $\boldsymbol{\kappa} = \kappa \mathbf{I}$, here considering isotropic diffusion, by the temperature gradient $\nabla \theta$. 

\subsection{Weak enforcement of boundary and interface conditions}\label{subsectionBoundary}
Boundary and interface conditions are imposed weakly via Nitsche's formulation, see \cite{1971Nitsche}. To enforce prescribed displacements on Dirichlet boundaries, the static equilibrium in Eq.~\eqref{Eq_RLinU} is augmented with:
\begin{equation}
\mathcal{R}^{\mathbf{u}}_{D} = \sum_I \Bigg[\ 
 - \int_{\Gamma^{I}_{D}} \delta \mathbf{u} \cdot \left( \boldsymbol{\sigma} \cdot \mathbf{n}_{\Gamma} \right) d\Gamma
 + \int_{\Gamma^{I}_{D}} \delta \left( \boldsymbol{\sigma} \cdot \mathbf{n}_{\Gamma} \right) \cdot \left( \mathbf{u} - \mathbf{u}_D \right)\ d\Gamma
 + \int_{\Gamma^{I}_{D}} \gamma^{\mathbf{u}}_{N}\ \delta \mathbf{u} \cdot \left( \mathbf{u} - \mathbf{u}_D \right)\ d\Gamma\ \Bigg],
\label{Eq_RNitscheU}
\end{equation}
where a non-symmetric Nitsche formulation is considered, see for example \cite{2012Burman} and \cite{2016SchillingerEtAlb}, and $\mathbf{u}_D$ is the displacement imposed on the Dirichlet boundary $\Gamma^{I}_{D}$. The parameter $\gamma^{\mathbf{u}}_{N}$ is chosen to achieve a desired accuracy in satisfying the boundary conditions and is a multiple of the ratio $E/h$, where $E$ is the Young's modulus of the considered material and $h$ is the edge length of the intersected background element.

The same formulation is used to impose temperature on Dirichlet boundaries by augmenting the equilibrium in Eq.~\eqref{Eq_RLinTheta} with:
\begin{equation}
\mathcal{R}^{\theta}_{D} = \sum_I \Bigg[\ 
 - \int_{\Gamma^{I}_{D}} \delta \theta\ \left( \mathbf{q} \cdot \mathbf{n}_{\Gamma} \right) d\Gamma
 + \int_{\Gamma^{I}_{D}} \delta \left( \mathbf{q} \cdot \mathbf{n}_{\Gamma} \right) \left( \theta - \theta_D \right)\ d\Gamma
 + \int_{\Gamma^{I}_{D}} \gamma^{\theta}_{N}\ \delta \theta\ \left( \theta - \theta_D \right)\ d\Gamma\ \Bigg],
\label{Eq_RNitscheTheta}
\end{equation}
where a non-symmetric Nitsche formulation is again considered and $\theta_D$ is the imposed temperature on the Dirichlet boundary $\Gamma^{I}_{D}$. The parameter $\gamma^{\theta}_{N}$ is defined similarly to $\gamma^{\mathbf{u}}_{N}$ as a multiple of the ratio $\kappa/h$, where $\kappa$ is the conductivity of the considered material.

Interface conditions between the displacements of materials $I$ and $J$ on $\Gamma^{IJ} = \partial \Omega^{I} \cap \partial \Omega^{J}$ are imposed similarly using the non-symmetric formulation: 
\begin{equation}
\mathcal{R}^{\mathbf{u}}_{Itf} = \sum_{\substack{I,J\\ I\neq J}} \Bigg[\
 - \int_{\Gamma^{IJ}} \delta \mathbf{u} \cdot \{ \boldsymbol{\sigma}\cdot \mathbf{n}_{\Gamma} \}\ d\Gamma
 + \int_{\Gamma^{IJ}} \delta \left( \boldsymbol{\sigma} \cdot \mathbf{n}_{\Gamma} \right) \cdot \llbracket \mathbf{u}\rrbracket\ d\Gamma
 + \int_{\Gamma^{IJ}} \gamma^{\mathbf{u}}_{Itf}\ \delta \mathbf{u} \cdot \llbracket \mathbf{u} \rrbracket\ d\Gamma\ \Bigg],
\label{Eq_RInterfaceU}
\end{equation}
where the jump operator $\llbracket \bullet \rrbracket$ computes the difference in the considered quantity between material domains $I$ and $J$ as $\llbracket \bullet \rrbracket = \bullet^{I} - \bullet^{J}$. The mean operator $\{ \bullet \}$ computes a weighted sum of the considered quantity over the materials $I$ and $J$ as $\{ \bullet \} = w^{I}\ \bullet^{I} + w^{J}\ \bullet^{J}$. The weights, following \cite{2009DolbowHarari} and \cite{2012AnnavarapuEtAl}, are defined as:
 \begin{equation}
w^{I} = \frac{\mbox{meas}(\Omega^{I})/E^{I}}{\mbox{meas}(\Omega^{I})/E^{I} + \mbox{meas}(\Omega^{J})/E^{J}},
\label{Eq_weightItf_1}
\end{equation}
and:
\begin{equation}
w^{J} = \frac{\mbox{meas}(\Omega^{J})/E^{J}}{\mbox{meas}(\Omega^{I})/E^{I} + \mbox{meas}(\Omega^{J})/E^{J}},
\label{Eq_weightItf_2}
\end{equation} 
where $\mbox{meas}(\Omega^{I})$ and $\mbox{meas}(\Omega^{J})$ are the volume or the surface area of materials $I$ and $J$ within the intersected element in two and three dimensions respectively. The properties $E^{I}$ and $E^{J}$ are the Young's moduli of materials $I$ and $J$. The penalty parameter $\gamma^{\mathbf{u}}_{Itf}$ is evaluated as: 
\begin{equation}
\gamma^{\mathbf{u}}_{Itf} = \frac{2\, E^{I}\, \mbox{meas}(\Gamma^{IJ})}{\mbox{meas}(\Omega^{I})/E^{I} + \mbox{meas}(\Omega^{J})/E^{J}},
\label{Eq_gammaItf}
\end{equation}
where $\mbox{meas}(\Gamma^{IJ})$ is a the surface area or the length of the interface within the intersected element in two or three dimensions respectively.

Interface conditions on the temperature between material domains $I$ and $J$ are imposed as follows:
\begin{equation}
\mathcal{R}^{\theta}_{Itf} = \sum_{\substack{I,J\\ I\neq J}} \Bigg[\
 - \int_{\Gamma^{IJ}} \delta \theta\ \{ \mathbf{q} \cdot \mathbf{n}_{\Gamma} \}\ d\Gamma
 + \int_{\Gamma^{IJ}} \delta \left( \mathbf{q} \cdot \mathbf{n}_{\Gamma} \right) \llbracket \theta \rrbracket\ d\Gamma
 + \int_{\Gamma^{IJ}} \gamma^{\theta}_{Itf}\ \delta \theta\ \llbracket \theta \rrbracket\ d\Gamma\ \Bigg],
\label{Eq_RInterfaceTheta}
\end{equation}
where the weights $w^{I}$, $w^{J}$ and the penalty parameter $\gamma^{\theta}_{Itf}$ are evaluated similarly to the linear elastic case as described in Eqs.~\eqref{Eq_weightItf_1}, \eqref{Eq_weightItf_2}, and \eqref{Eq_gammaItf}, but by substituting the conductivity $\kappa$ for the Young's modulus $E$ of the considered material.

\subsection{Face-oriented ghost stabilization}\label{subsectionStabilization}
Using enriched finite element techniques, numerical instabilities arise when either the contributions of one or more of the basis functions approximating the physics variable field to the residual vanish and/or these contributions become linearly dependent. These issues typically arise when the level set field intersects elements such that small material subdomains emerge. This results in an ill-conditioning of the equation system and inaccurate prediction of solution field gradients along the interface. Different techniques are available to mitigate this issue; for example ghost stabilization (see \cite{2010Burman}), basis function removal (see \cite{2018ElfversonEtAl} for the specific case of B-spline interpolation), basis function aggregation (see \cite{2022BadiaEtAl}). In this work, a generalized version of the face-oriented ghost stabilization, proposed by \cite{2014BurmanHansbo}, is used to fit the proposed enrichment strategy based on the basis functions supports as described hereunder.

\begin{figure*}[!h]\center
\includegraphics[scale=1]{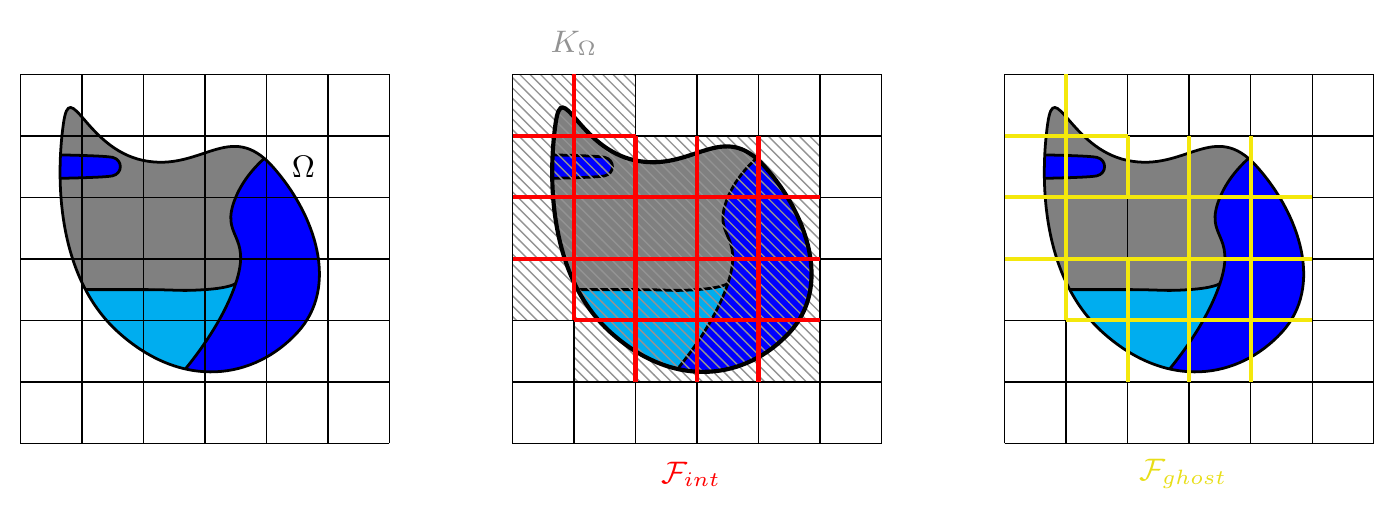}
\caption{Set of interior facets $\mathcal{F}_{int}$ and set of ghost facets $\mathcal{F}_{ghost}$ used for the face-oriented ghost stabilization for a three-material problem.}
\label{figGhost0}
\end{figure*}

The domain $\Omega$, made of the union of all material domains $\Omega^{I}$, is immersed in a background mesh. The set of all background elements in the mesh is denoted $\mathcal{K}$ and $K_{\Omega}$ is the subset of background elements that have a non-empty intersection with $\Omega$:
\begin{equation}
K_{\Omega} \coloneqq \left\{ K \in \mathcal{K} : K \cap \Omega \neq \emptyset \right\},
\label{Eq_ghost_K}
\end{equation}
The set of interior facets of $K_{\Omega}$ is denoted $\mathcal{F}_{int}$. Each interior facet $F \in \mathcal{F}_{int}$ is shared between two elements $\Omega_{F}^{+}$ and $\Omega_{F}^{-}$ of $K_{\Omega}$. Finally $\tilde{\Gamma}$ is defined as the union of all material interfaces and geometric boundaries. The set of ghost facets is then taken to be:
\begin{equation}
\mathcal{F}_{ghost} \coloneqq \Big\{ F \in \mathcal{F}_{int} : \Omega_{F}^{+} \cap \tilde{\Gamma} \neq \emptyset,
\mbox{or} \ \Omega_{F}^{-} \cap \tilde{\Gamma} \neq \emptyset \Big\}.
\label{Eq_ghost_F}
\end{equation}
The set of interior facets $\mathcal{F}_{int}$ and the set of ghost facets $\mathcal{F}_{ghost}$ are illustrated on a three-material problem in Fig. \ref{figGhost0}. A three-material domain $\Omega$ is immersed in a background mesh. The set $K_{\Omega}$ of background elements with non-empty intersection with $\Omega$ are hatched in grey. The facets that lie within $K_{\Omega}$ form the set $\mathcal{F}_{int}$ and are drawn in red. Finally, the set of ghost facets is shown in yellow.

\begin{figure*}[!h]\center
\includegraphics[scale=1]{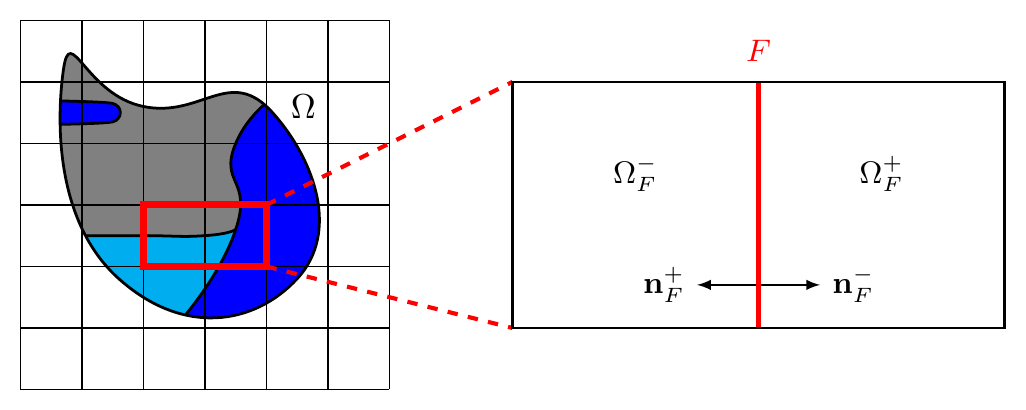}
\caption{Description of ghost facet $F$ shared between two adjacent background elements $\Omega_{F}^{+}$ and $\Omega_F^{-}$ for the face-oriented ghost stabilization for a three-material problem.}
\label{figGhost1}
\end{figure*}

Consider ghost facet $F \in \mathcal{F}_{ghost}$ shared between two adjacent background elements $\Omega_{F}^{+}$ and $\Omega_F^{-}$ as illustrated in Fig. \ref{figGhost1}. The outward facing normals to $\Omega_{F}^{+}$ and $\Omega_{F}^{-}$ along $F$ are defined as $\mathbf{n}_{F}^{+}$ and $\mathbf{n}_{F}^{-}$, such that $\mathbf{n}_{F} = \mathbf{n}_{F}^{+} = - \mathbf{n}_{F}^{-}$. The material layout subdivides the element $\Omega_{F}^{+}$ into $N_{F}^{+}$ connected subdomains $\{\Omega_{F,i}^{+}\}_{i=1}^{N_{F}^{+}}$ and the element $\Omega_{F}^{-}$ into $N_{F}^{-}$ connected subdomains $\{ \Omega_{F}^{-} \}_{j=1}^{N_{F}^{-}}$ such that each subdomain is occupied by one and only one material. 

The material index associated to $\Omega_{F,i}^{+}$ for each $i=1,\dots,N_{F}^{+}$ is denoted $\mathcal{M}_{F,i}^{+} = \mathcal{M}(\Omega_{F,i}^{+})$ and the material index associated to $\Omega_{F,j}^{-}$ for each $j=1,\dots,N_{F}^{-}$ is $\mathcal{M}_{F,j}^{-} = \mathcal{M}(\Omega_{F,j}^{-})$. Finally, the polynomial extension of the field $\mathbf{u}|_{\Omega_{F,i}^{+}}$ to all of $\mathbb{R}^{d}$ for each $i=1,\dots,N_F^{+}$ is defined as $\mathbf{u}_{F,i}^{+}$ and the polynomial extension of the field $\mathbf{u}|_{\Omega_{F,j}^{-}}$ to all of $\mathbb{R}^{d}$ for each $j=1,\dots,N_F^{-}$ as $\mathbf{u}_{F,j}^{-}$. The ghost stabilization $\mathcal{G}^{\mathbf{u}}_{F}$ for facet $F$, penalizing the jumps in the displacement gradients across the facet, is then taken to be:  
\begin{equation}
\mathcal{G}^{\mathbf{u}}_{F}(\mathbf{u}, \delta \mathbf{u}) =
\sum_{i=1}^{N_{F}^{+}} 
\sum_{j \in J_{F,i}}
\Bigg[ \sum_{k=1}^{p} \int_{F}  
\gamma_{G}^{\mathbf{u}}\ h^{\tilde{k}} 
\Big\llbracket \partial^{k}_{n}\, \delta \mathbf{u} \Big\rrbracket 
\cdot
\Big\llbracket \partial^{k}_{n}\, \mathbf{u} \Big\rrbracket\, 
d\Gamma \Bigg],
\label{Eq_RGhostU1}
\end{equation}
where the set $J_{F,i}$ is defined so that:
\begin{equation}
J_{F,i} \coloneqq \Big\{ j \in \{ 1,\dots, N_{F}^{-}  \} : \mathcal{M}_{F,i}^{+}=\mathcal{M}_{F,j}^{-}\neq0,
\mbox{and}\ | \partial \Omega_{F,i}^{+}\, \cap\, \partial \Omega_{F,j}^{-} | \neq 0 \Big\}.
\label{Eq_JSet}
\end{equation}
The parameter $\tilde{k} = {2(k-1)+1}$ and $\llbracket \bullet \rrbracket$ is a jump operator such that:
\begin{equation}
\Big\llbracket \partial^{k}_{n}\, \delta \mathbf{u} \Big\rrbracket = \left( \partial^{k}_{n} \delta \mathbf{u}_{F,i}^{+} - \partial^{k}_{n} \delta \mathbf{u}_{F,j}^{-} \right),
\end{equation}
and:
\begin{equation}
\Big\llbracket \partial^{k}_{n}\, \mathbf{u} \Big\rrbracket = \left( \partial^{k}_{n} \mathbf{u}_{F,i}^{+} - \partial^{k}_{n} \mathbf{u}_{F,j}^{-} \right).
\end{equation}
The operator $\partial_n^k(\bullet)$ is the $k^{th}$ order normal derivative operator and $\partial_n(\bullet) = \nabla(\bullet) \cdot \mathbf{n}_F$ with $\nabla(\bullet)$ the spatial derivative. The parameter $p$ is the degree of the considered approximation. It should be noted that only the $p^{th}$ contribution is non-zero for maximally smooth splines, i.e., splines without repeating knots. The ghost penalty parameter $\gamma_{G}^{\mathbf{u}}$ is defined as a multiple of the Young's modulus $E$ of the considered material. 

The formulation in Eq.~\eqref{Eq_JSet} requires that:
$$\mathcal{M}_{F,i}^{+}=\mathcal{M}_{F,j}^{-}\neq0,$$
assuming that the material index for void is zero, so that the ghost stabilization is only applied between $\mathbf{u}_{F,i}^{+}$ and $\mathbf{u}_{F,j}^{-}$ when $\Omega_{F,i}^{+}$ and $\Omega_{F,j}^{-}$ are occupied by the same non-void material. Additionally, the formulation requires that:
$$| \partial \Omega_{F,i}^{+}\, \cap\, \partial \Omega_{F,j}^{-} | \neq 0,$$
and the ghost stabilization is only applied between $\mathbf{u}_{F,i}^{+}$ and $\mathbf{u}_{F,j}^{-}$ when the boundaries of $\Omega_{F,i}^{+}$ and $\Omega_{F,j}^{-}$, $\partial \Omega_{F,i}^{+}$ and $\partial \Omega_{F,j}^{-}$ respectively, meet along a portion of the facet $F$ with a non-zero measure, e.g., the boundaries meet along more than a point in two dimensions and along more than a line in three dimensions. 

\begin{figure*}[!h]\center
\includegraphics[scale=1]{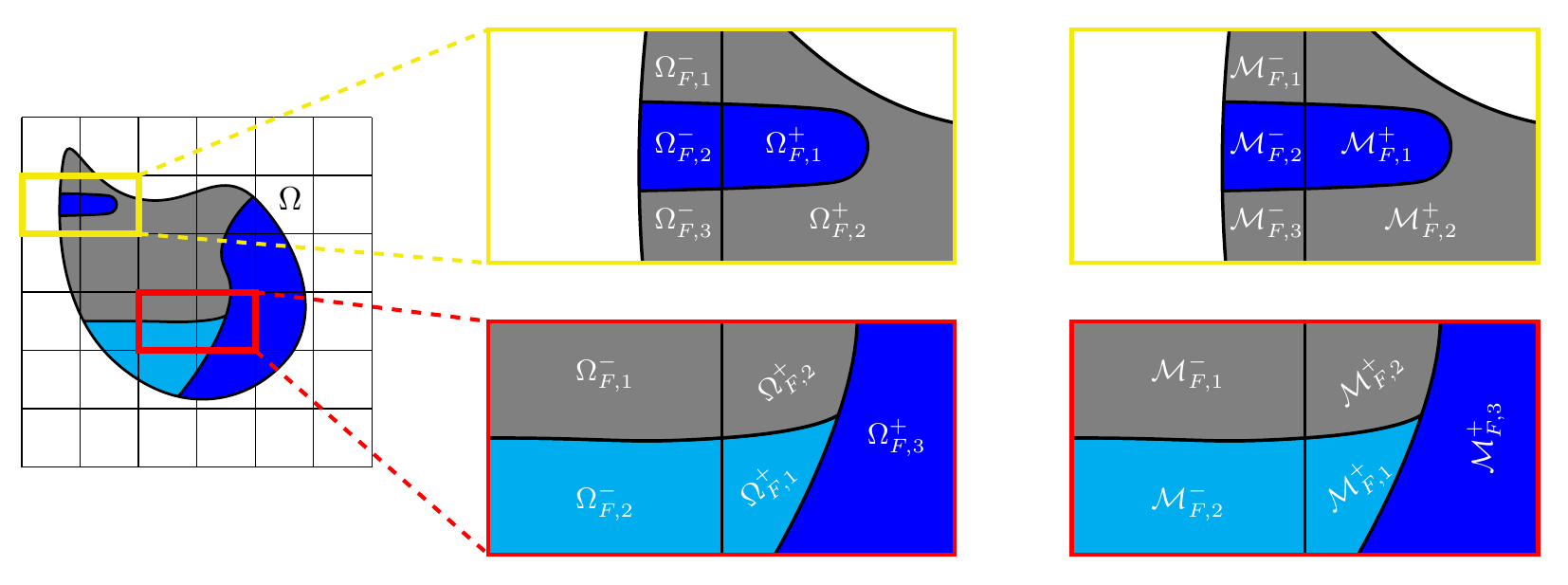}
\caption{Subdivision of $\Omega_{F}^{+}$ and $\Omega_{F}^{+}$ into connected subdomains $\Omega_{F,i}^{+},\ i=1, \dots, N_F^{+}$ and $\Omega_{F,j}^{-},\ j=1, \dots, N_F^{-}$ for the face-oriented ghost stabilization for a three-material problem.}
\label{figGhost2}
\end{figure*}

The subdivision of $\Omega_{F}^{+}$ and $\Omega_{F}^{-}$ into connected subdomains $\Omega_{F,i}^{+}$ and $\Omega_{F,j}^{-}$ and the associated material indices $\mathcal{M}_{F,i}^{+}$ and $\mathcal{M}_{F,j}^{-}$ is shown in Fig. \ref{figGhost2} for two different material configurations. In the first case marked by a red box, the background element $\Omega_{F}^{-}$ is occupied by two connected material subdomains $\Omega_{F,1}^{-}$ and $\Omega_{F,2}^{-}$, while the background element $\Omega_{F}^{+}$ is divided into three connected material subdomains $\Omega_{F,1}^{+}$, $\Omega_{F,2}^{+}$, and $\Omega_{F,3}^{+}$. As $\Omega_{F,1}^{-}$ and $\Omega_{F,2}^{+}$ have the same material index $\mathcal{M}_{F,1}^{-} = \mathcal{M}_{F,2}^{+}$ and meet along the facet, the jump in the associated field gradients across is penalized across the facet. The same holds for $\Omega_{F,2}^{-}$ and $\Omega_{F,1}^{+}$. 

In the second case marked by a yellow box, a different situation arises. While the background element $\Omega_{F}^{-}$ is filled with only two different materials, the grey material lies within two connected subdomains within the element. The background element $\Omega_{F}^{-}$ is thus divided into three connected material subdomains $\Omega_{F,1}^{-}$, $\Omega_{F,2}^{-}$, and $\Omega_{F,3}^{-}$, while the background element $\Omega_{F}^{+}$ is divided into two connected material subdomains $\Omega_{F,1}^{+}$ and $\Omega_{F,2}^{+}$. As $\Omega_{F,1}^{-}$ and $\Omega_{F,2}^{+}$ have the same material index, $\mathcal{M}_{F,1}^{-} = \mathcal{M}_{F,2}^{+}$, and meet along the facet, the jump in the associated field gradients is penalized across the facet. The same holds for the pairs $\Omega_{F,3}^{-}$ and $\Omega_{F,2}^{+}$, and $\Omega_{F,2}^{-}$ and $\Omega_{F,1}^{+}$.

The full contribution of the  ghost stabilization to the residual equations is attained by summing over all ghost facets:
\begin{equation}
\mathcal{R}^{\mathbf{u}}_{Ghost} = 
\sum_{F \in \mathcal{F}_{ghost}} \mathcal{G}^{\mathbf{u}}_{F}(\mathbf{u}, \delta \mathbf{u}).
\label{Eq_RGhostU2}
\end{equation}

The ghost penalization for the temperature field is defined similarly as:
\begin{equation}
\mathcal{R}^{\theta}_{Ghost} = 
\sum_{F \in \mathcal{F}_{ghost}} 
\mathcal{G}^{\theta}_{F}(\theta, \delta \theta),
\label{Eq_RGhostTheta}
\end{equation}
where the ghost stabilization $\mathcal{G}_{F}^{\theta}$ for facet $F$ is:
\begin{equation}
\mathcal{G}^{\theta}_{F}(\theta, \delta \theta) = 
\sum_{i=1}^{N_{F}^{+}} 
\sum_{j \in J_{F,i}} 
\Bigg[ \sum_{k=1}^{p} \int_{F}  
\gamma_{G}^{\theta}\ h^{\tilde{k}} 
\Big\llbracket \partial^{k}_{n}\, \delta \theta \Big\rrbracket 
\cdot
\Big\llbracket \partial^{k}_{n}\, \theta \Big\rrbracket\, 
d\Gamma\ \Bigg],
\end{equation}
with the ghost penalty parameter, $\gamma_{G}^{\theta}$, is defined as a multiple of the conductivity $\kappa$ for the considered material.

\subsection{Numerical integration}\label{subsectionIntegration}
Working with immersed boundary techniques and using the Heaviside enrichment, the weak form of the governing equations is integrated separately on each material subdomain. Elements occupied by more than one material are decomposed into conforming integration subdomains. In two (three) dimensions, a quadrangle (a hexahedron) is subdivided into a triangular (tetrahedral) integration mesh that conforms to the material interfaces. Gauss quadrature rules are used on the generated integration elements.

\begin{figure*}[!h]\center
\includegraphics[scale=1]{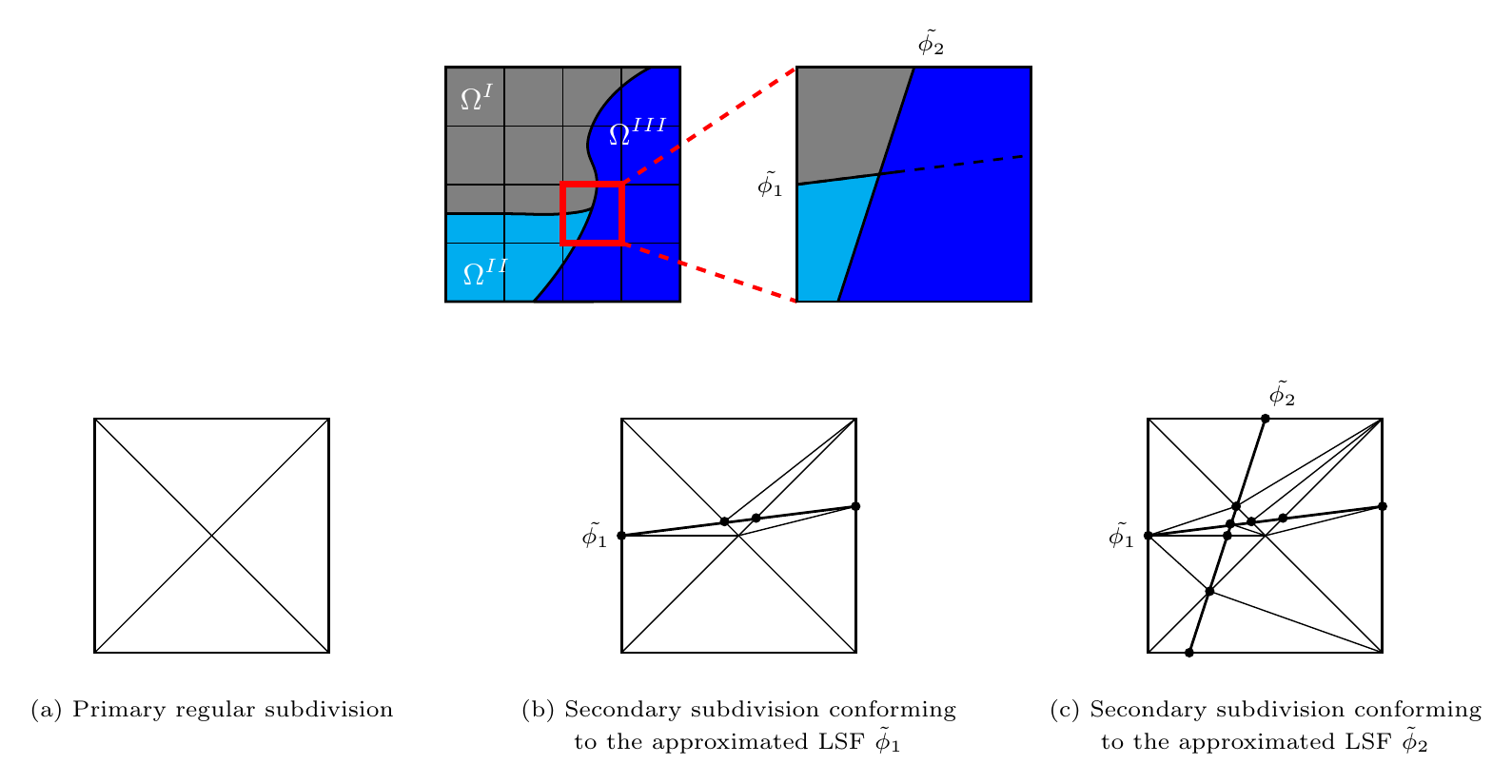}
\caption{Generation of conforming numerical integration mesh using a primary and a secondary subdivision for a three-material problem.}
\label{figNumericalIntegration}
\end{figure*}

The subdivision strategy is illustrated with a two dimensional three-material problem in Fig. \ref{figNumericalIntegration}. First, to increase the accuracy of the interface detection, a primary subdivision is performed, and the background element is divided in four subtriangles. Then, a secondary subdivision is performed to construct a triangular mesh that conforms to the interfaces created by the LSFs. The LSFs are linearly interpolated along the element edges using $\tilde{\phi}_1$ and $\tilde{\phi}_2$. The intersections between the approximated LSFs and the element edges are determined, see black circles on the figure, and triangular integration elements are created. It should be noted that additional refinement of the background mesh can be carried out before constructing the integration mesh from the LSFs to achieve reduced geometric error, see Subsection \ref{Ex4}.

\section{Numerical examples}\label{sectionExamples}
In this section, the versatility of the proposed XIGA approach and its ability to systematically and accurately address multi-material problems is demonstrated. First, the accuracy of the evaluated responses is investigated  in Subsection \ref{Ex1}. The influence of the minimum size of material integration subdomains within a basis function and of the choice of the ghost penalty parameter on the conditioning of the system of equations is studied. Subsection \ref{Ex2} focuses on the robustness of the approach with respect to the creation of small material integration subdomains by varying the location of the analysis domain within a fixed background mesh. In Subsection \ref{Ex3}, the approach is used to solve \textit{N}-phase junction problems, and the convergence rates attained with \textit{h}-refinement considering linear, quadratic, and cubic B-spline basis functions are investigated. The ability of the XIGA approach to handle non-planar interface configurations is assessed in Subsection \ref{Ex4}. Finally, a \textit{N}-material problem is tackled in Subsection \ref{Ex5}.

In the following examples, the performances of the proposed XIGA approach are studied and characterized using three criteria: the system condition number, the relative $L^2$ error norm and the relative $H^1$ error semi-norm. 
Ill-conditioning can affect the convergence of the solver of the system of equations. The condition number is used to assess the conditioning of the system of equations and is evaluated as:  
\begin{equation}
\mbox{cond}(\mathbf{A}) = || \mathbf{A}^{-1} || \cdot || \mathbf{A} ||,
\label{Eq_Cond}
\end{equation}
where $\mathbf{A}$ is a matrix describing the system of equations to solve and $|| \bullet ||$ is the Frobenius norm.

For a generic state field $\mathbf{a}$, the relative $L^2$ error norm is defined as:
\begin{equation}\displaystyle
L^2 = 
\sqrt{ 
\frac{\int_{\Omega} \left( \mathbf{a}^h - \mathbf{a} \right)^T\, \left( \mathbf{a}^h - \mathbf{a} \right)\ d \Omega}
{\int_{\Omega} \mathbf{a}^T\, \mathbf{a}\ d \Omega}},
\label{Eq_L2}
\end{equation}
where $\mathbf{a}$ is a reference solution, here chosen as either an analytical solution or a numerical solution computed on a sufficiently refined mesh, and $\mathbf{a}^h$ is the numerical solution evaluated with the XIGA approach.

For a generic state field $\mathbf{a}$, the relative $H^1$ error semi-norm is evaluated as:
\begin{equation}\displaystyle
H^1 = 
\sqrt{ 
\frac{\int_{\Omega} \left( \nabla\mathbf{a}^h - \nabla \mathbf{a} \right)^T\, \left(\nabla\mathbf{a}^h - \nabla \mathbf{a} \right)\ d \Omega}
{\int_{\Omega} \nabla \mathbf{a}^T\, \nabla \mathbf{a}\ d \Omega}},
\label{Eq_H1}
\end{equation}
where $\nabla \mathbf{a}$ is the gradient of a reference solution, here chosen as either an analytical solution or a numerical solution computed on a sufficiently refined mesh, and $\nabla \mathbf{a}^h$ is the gradient of the numerical solution evaluated with the XIGA approach.

In all the following examples, the set of discretized governing equations is integrated using Gauss quadrature rules on each integration subelement depending on the order of the basis functions. In two dimensions, $2{\times}2$-, $3{\times}3$-, and $4{\times}4$-point integration rules are used for quadrangular integration elements and 7-, 12-, or 25-point integration rules are used for triangular integration elements for linear, quadratic, and cubic basis functions respectively. In three dimensions, $2{\times}2{\times}2$-, $3{\times}3{\times}3$-, and $4{\times}4{\times}4$-point integration rules are used for hexahedral integration elements and 11-, 35-, and 56-point integration rules are used for tetrahedral integration elements for linear, quadratic, and cubic basis functions respectively. The systems of discretized governing equations are built using an implementation of the XIGA approach within an in-house fully parallelized C++ code, and are solved by the direct solver PARDISO for the 2D problems (see \cite{2018KourounisEtAl}), and by a GMRES algorithm for 3D problems, preconditioned by an algebraic multi-grid solver (see \cite{2006GeeEtAl}). In some of the numerical examples, the systems of equations are poorly conditioned and condition numbers exceeding $10^{25}$ are observed. Despite these large condition numbers, the linear solve converged for all presented results.

\subsection{Stability study with respect to material subdomains size and ghost penalty parameter}\label{Ex1}
In this subsection, the accuracy of the XIGA approach is demonstrated by showing that low errors with respect to the analytical solution can be achieved. The conditioning of the system of equations with respect to the size of the created material integration subdomains and the value of the ghost penalty parameter is studied.

\begin{figure*}[!h]\center
\includegraphics[scale=1]{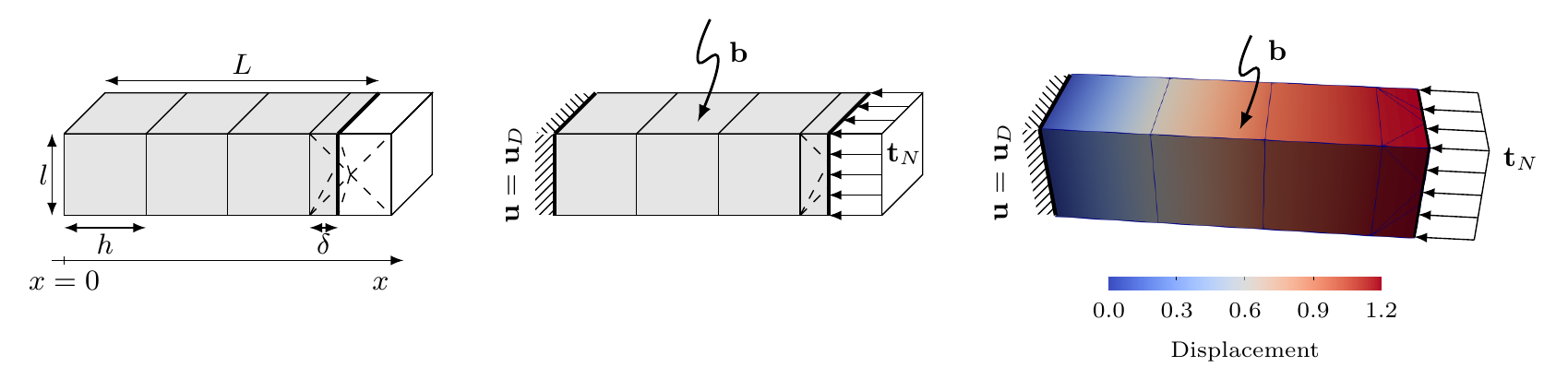}
\caption{Single material linear elastic bar with varying sliver size and ghost penalty parameter: problem setup, boundary conditions and solution for quadratic state solution using quadratic B-splines.}
\label{Ex1_fig0}
\end{figure*}

For this purpose, a single material bar is considered in three dimensions, with dimensions $L = 3 + \delta\, \mbox{m}$, $l = 1\, \mbox{m}$, and a cross-section area $A = l^2 = 1\, \mbox{m}^2$. We consider a linear elastic problem defined on this geometry. The set up and boundary conditions are illustrated in Fig. \ref{Ex1_fig0}. Considering linear elasticity, the material Young's modulus is set to $E = 10.0\ \mbox{N/m}^2$ and the Poisson ration is set to $\nu = 0.0$ to avoid any three-dimensional effect. The left side of the bar is clamped, i.e., $\mathbf{u}_D = [u_{Dx}\ u_{Dy}\ u_{Dz}]^T = [0.0\ 0.0\ 0.0]^T\, \mbox{m}$, and the Nitsche's penalty parameter is set to $\gamma_N = 100.0$. Three loading scenarios are considered: (i) a traction at the bar tip $\mathbf{t}_N = [t_{Nx}\ t_{Ny}\ t_{Nz}]^T = [5.0\ 0.0\ 0.0]\, \mbox{N/m}^2$, (ii) a constant body load $\mathbf{b} = b_0 [1.0\ 0.0\ 0.0]\, \mbox{N/m}^3$ with $b_0 = 2.0$ within the material, and (iii) a linear body load $\mathbf{b} = b_0 [x\ 0.0\ 0.0]\, \mbox{N/m}^3$ with $b_0 = 2.0$ within the material. These three loading scenarios result in a one-dimensional linear elasticity problem and lead respectively to a linear, a quadratic, or a cubic displacement solution in $x$, that is the position in the horizontal direction:
\begin{equation}
\begin{array}{lll}
u_x(x) & = & \left\{
\begin{array}{ll}
\displaystyle u_{Dx} + \frac{t_{Nx}}{EA} x,& \mbox{for}\ t_{Nx} \neq 0.0,\, b_x = 0.0,\\[5pt]
\displaystyle u_{Dx} + \frac{b_0}{2 EA} \left( 2Lx - x^2 \right),& \mbox{for}\ t_{Nx} = 0.0,\, b_x = b_0, \\[5pt]
\displaystyle u_{Dx} + \frac{b_0}{6 EA}  \left( 3L^2x - x^3 \right),& \mbox{for}\ t_{Nx} = 0.0,\, b_x = b_0 x,
\end{array}
\right.\\
u_y & = & \quad 0,\\
u_z & = & \quad 0.
\end{array}
\label{Ex1_Eq1_analyticalSolution}
\end{equation}

Different locations of the end of the bar are considered generating different sliver sizes $\delta$ over the last element. The following sliver sizes are investigated: 
\begin{equation*}
\begin{split}
\delta = \big[ & 0.001\ 0.002\ 0.0035\ 0.005\ 0.007\ 0.01\ 0.015\ 0.025\ \\
& 0.04\ 0.06\ 0.08\ 0.1\ 0.15\ 0.25\ 0.4\ 0.6\ 0.8\ 0.9 \big] \times h,
\end{split}
\end{equation*}
where $h$ is the background element size and is set to $h = 1\, \mbox{m}$. Additionally, several values are considered for the ghost penalty parameter: 
$\gamma_G = \left[\, 0.0\ 10^{-9}\ 10^{-5}\ 10^{-4}\ 10^{-3}\ 10^{-2}\ 10^{-1}\ 1.0\, \right].$ The displacement fields are interpolated using linear, quadratic, and cubic B-spline basis functions. 

\begin{figure*}[!t]\center
\includegraphics[scale=0.75]{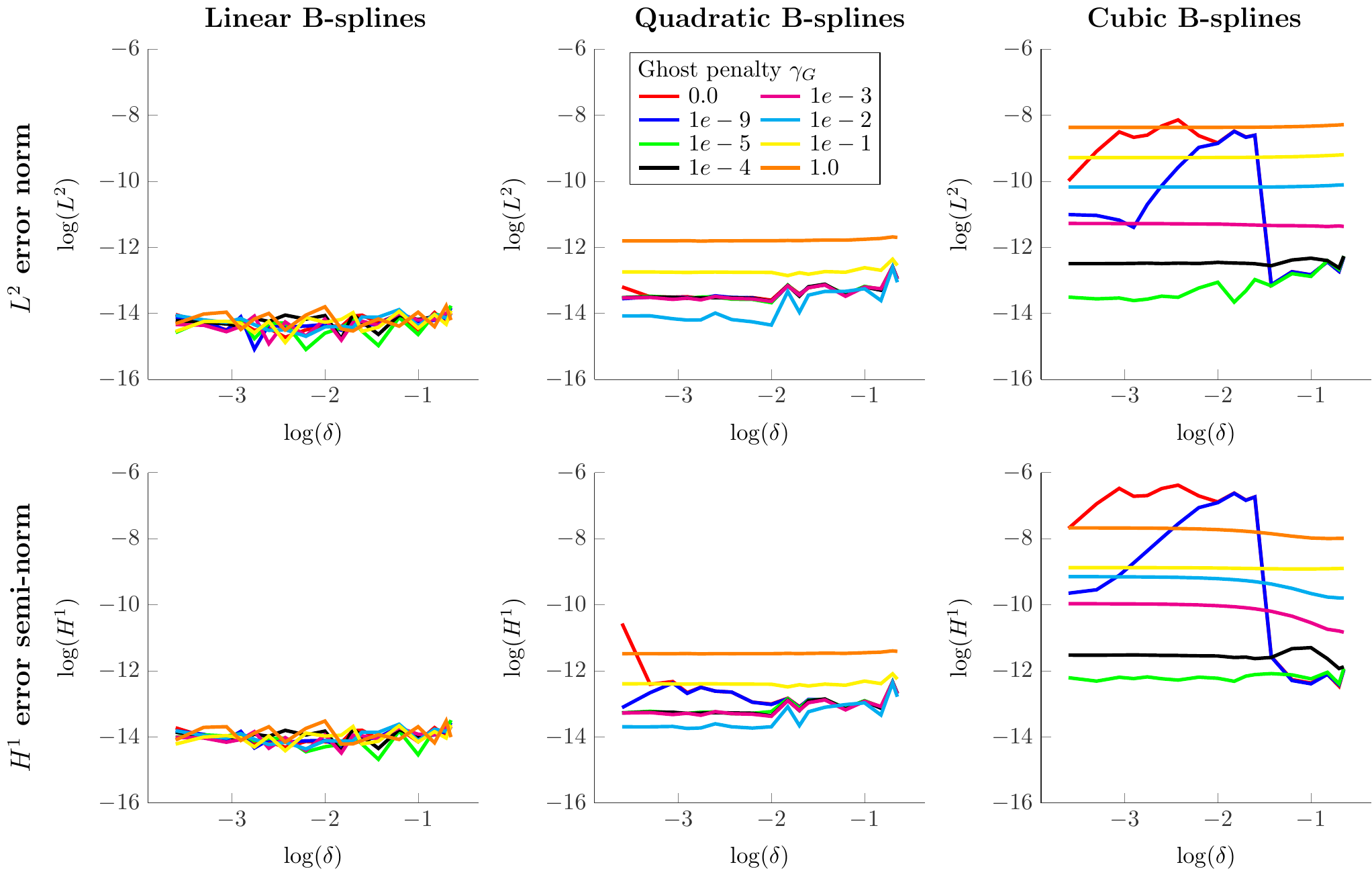}
\caption{Single material linear elastic bar undergoing a constant load at its tip, leading to a linear displacement solution.}
\label{Ex1_fig2}
\end{figure*}

\begin{figure*}[!h]\center
\includegraphics[scale=0.75]{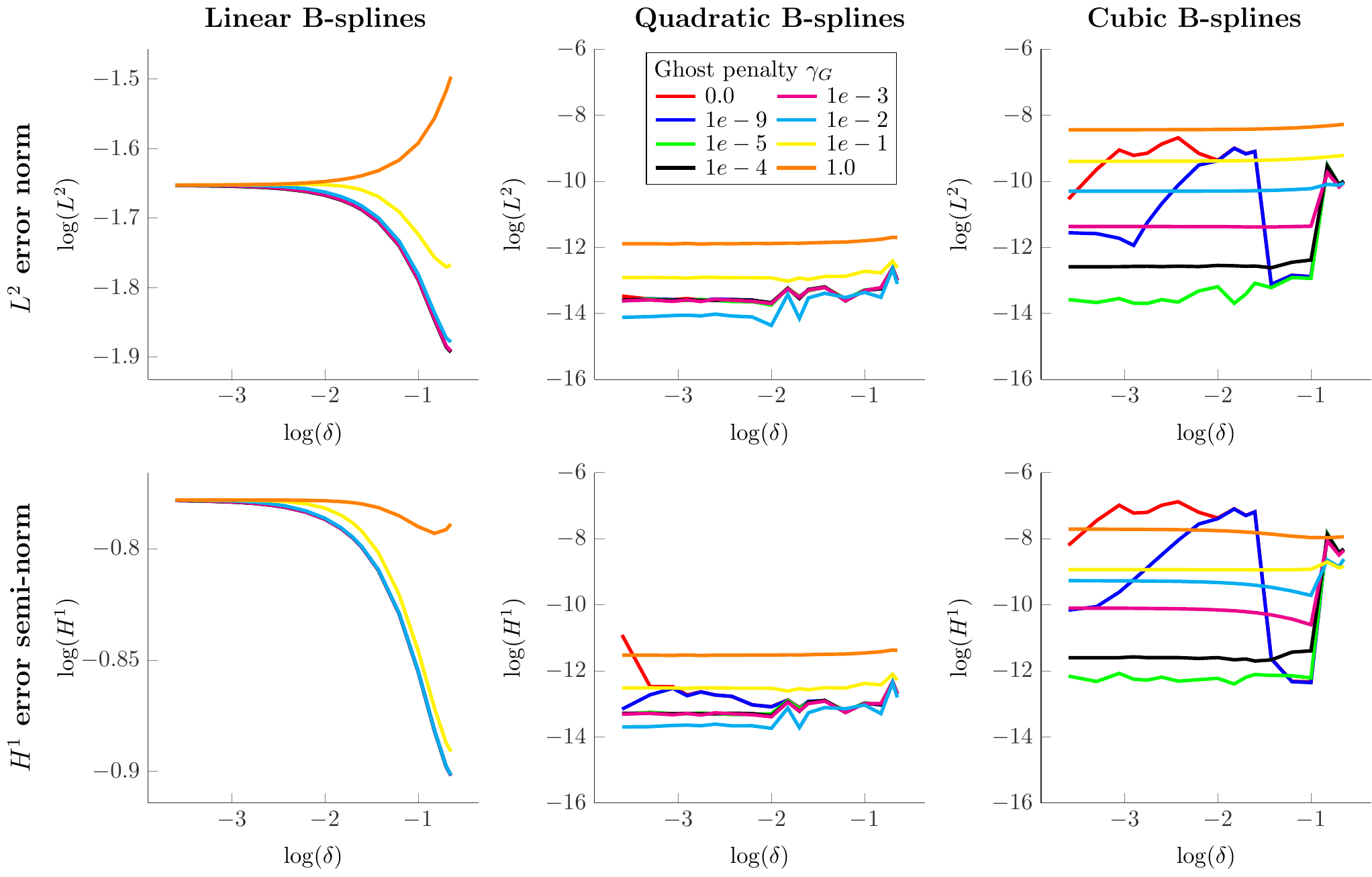}
\caption{Single material linear elastic bar undergoing a constant body load, leading to a quadratic displacement solution.}
\label{Ex1_fig3}
\end{figure*}

\begin{figure*}[!t]\center
\includegraphics[scale=0.75]{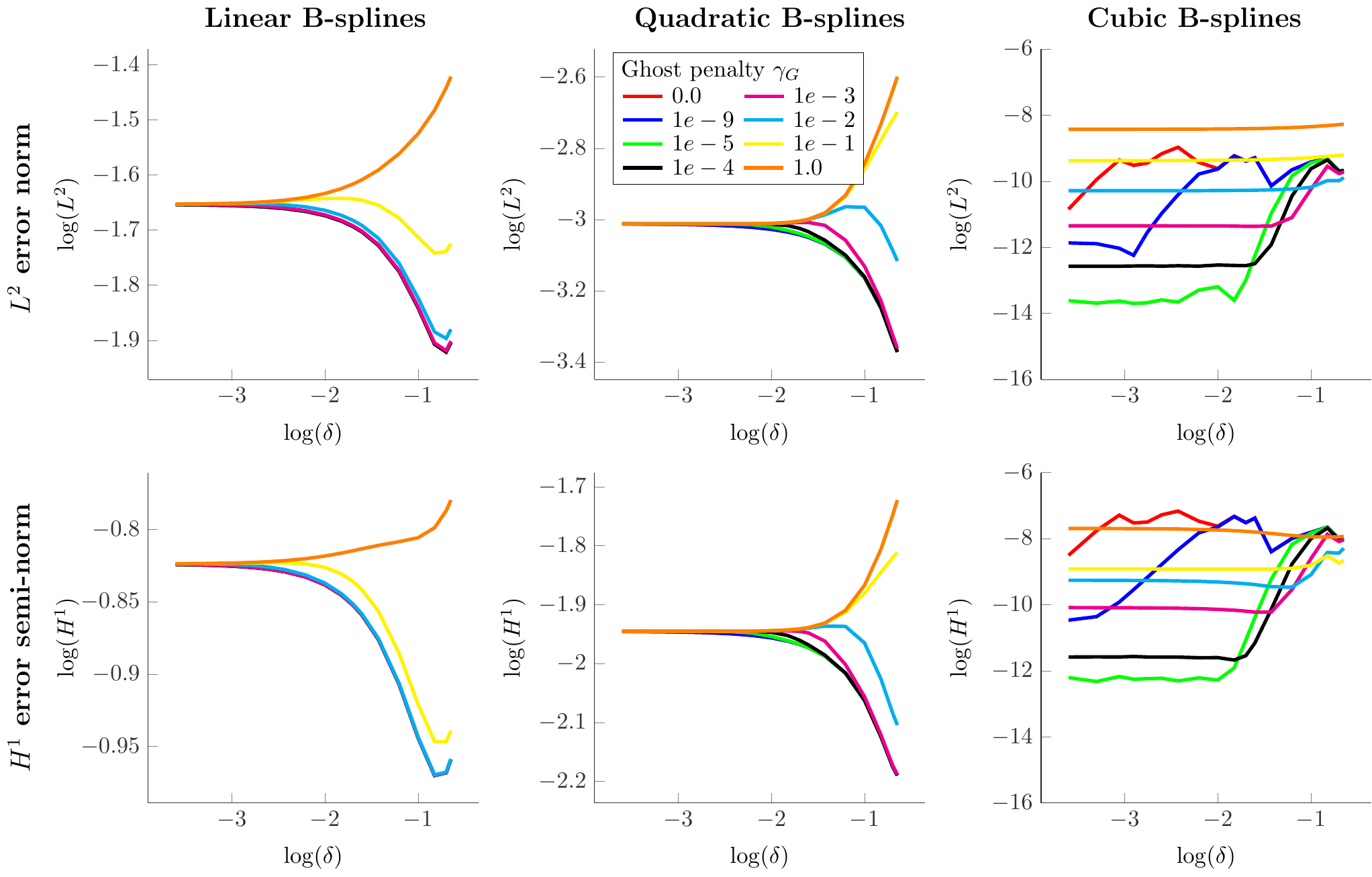}
\caption{Single material linear elastic bar undergoing a linear body load, leading to a cubic displacement solution.}
\label{Ex1_fig4}
\end{figure*}

\begin{figure*}[!h]\center
\includegraphics[scale=0.75]{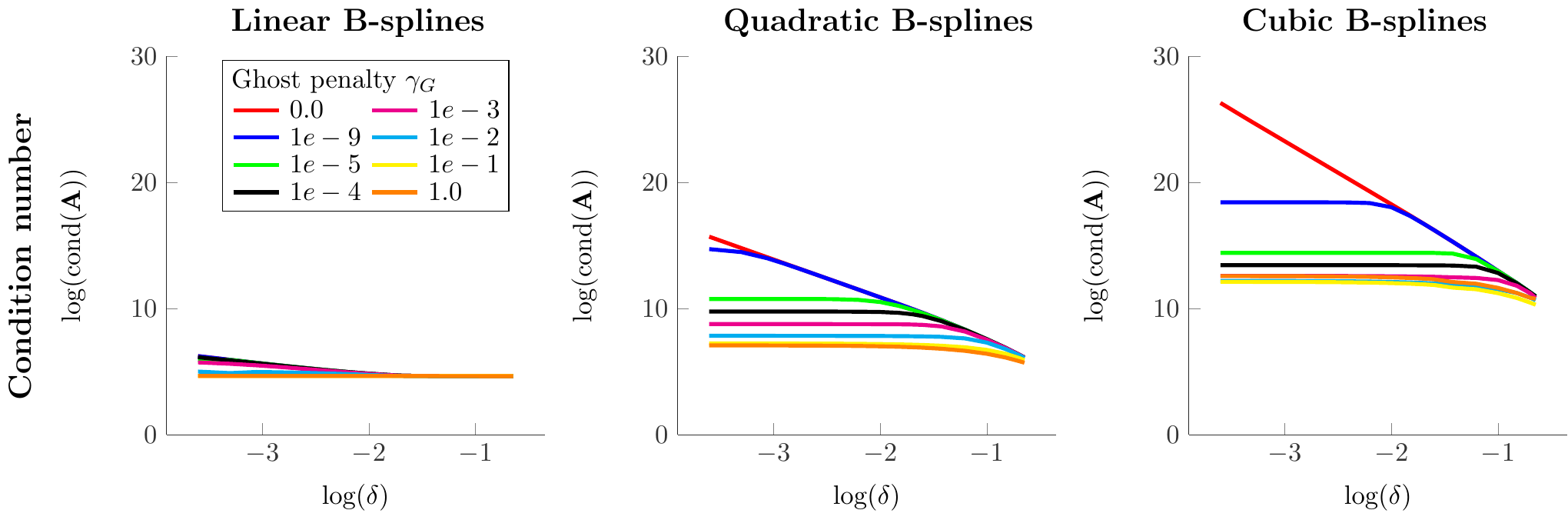}
\caption{Single material linear elastic bar undergoing a constant load at its tip, leading to a linear displacement solution along the bar.}
\label{Ex1_fig1}
\end{figure*}

The relative $L^2$ error norm and the relative $H^1$ error semi-norm are shown in Fig. \ref{Ex1_fig2} for the linear solution case, in Fig. \ref{Ex1_fig3} for the quadratic solution case, and in Fig. \ref{Ex1_fig4} for the cubic solution case. In these figures, each column presents a different interpolation order for the B-spline basis functions: linear, quadratic, and cubic. The rows display the relative $L^2$ error norm and the relative $H^1$ error semi-norm. Each ghost penalty value $\gamma_G$ is associated with a colored curve.
Focusing on Fig. \ref{Ex1_fig2}, the $L^2$ error norm and the $H^1$ error semi-norm values show that high accuracy can be achieved regardless of the B-spline order. This is expected as the exact solution is linear in $x$. In this case, the finite element solution is insensitive to the ghost penalty parameter choice or of the sliver size. Two additional observations are worth to be noted. A slight degradation of the $L^2$ and $H^1$ error occurs when the ghost penalty parameter is chosen too large, i.e., $\gamma_G \geq 1e-2$. In this case, the ghost stabilization acts as a coarsening operator and using high penalty values leads to errors similar to a coarsening of the mesh. When the ghost stabilization is turned off, i.e., $\gamma_G = 0$, or for small values of the penalty parameter, i.e.,$\gamma_G \leq 1e-5$, a slight effect of the sliver size is observed on the $L^2$ and $H^1$ errors, namely it increases as the sliver vanishes as expected. Low ghost penalty leads to poorly conditioned systems, as further observed with the condition number in Fig. \ref{Ex1_fig1}, and in turn to lower accuracy of the solution.

The results associated with the second and third loading cases leading to a quadratic and a cubic displacement solution are presented in Fig. \ref{Ex1_fig3} and \ref{Ex1_fig4} respectively. The results support the observations made for the linear solution case. It should be noted that for the second case, linear B-splines are not sufficient to accurately represent the exact solution that is quadratic in $x$. For the third case, both the linear and quadratic B-splines are not able to capture accurately the exact solution that is cubic in $x$. These results suggest that our XIGA approach is accurate and leads to low errors with respect to the analytical solution if the basis function order is sufficient to represent the analytical solution.

The condition number $\mbox{cond}(\mathbf{A})$ is the same for all loading scenarios. The condition numbers are displayed in Fig. \ref{Ex1_fig1} for different interface configurations and values of the ghost penalty parameter $\gamma_G$. Figure \ref{Ex1_fig1} shows an increase in the condition number when higher-order bases are used. This is due to the increased number of basis functions that are supported on intersected elements, see \cite{2017dePrenterEtAl}. Using ghost stabilization, the condition number can be significantly improved by using a penalty parameter value $\gamma_G \geq 1e-5$. This effect is particularly visible when using higher-order B-spline basis functions, i.e., quadratic and cubic. 

\subsection{Robustness study with respect to the creation of arbitrary integration subelements}\label{Ex2}
\begin{figure*}[!h]\center
\includegraphics[scale=1]{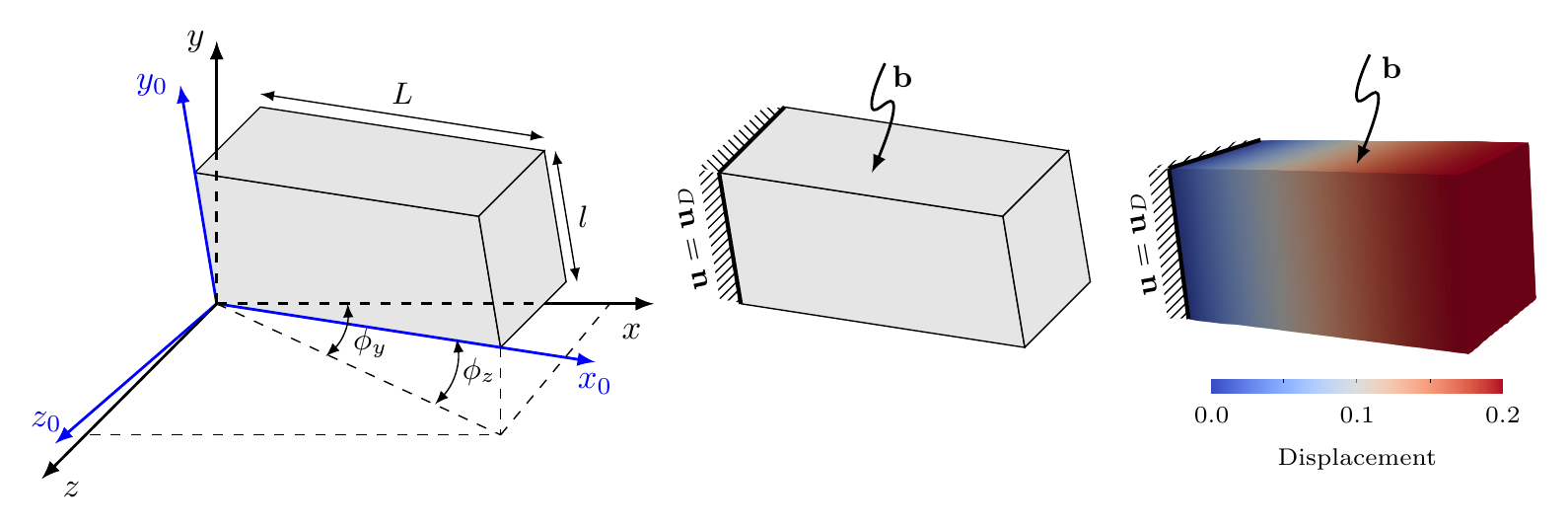}
\caption{An immersed single elastic material bar rotated in a fixed 3D background mesh: problem setup, boundary conditions and solution using quadratic B-spline with rotation angles $\phi_y = 20^o$, $\phi_z = 20^o$ and a mesh size $h=0.125\, \mbox{m}$.}
\label{Ex2_fig1}
\end{figure*}

To study the robustness of the method and the effectiveness of the stabilization with respect to different intersection configurations, a straight bar is immersed and rotated in a fixed background mesh. The problem is solved for linear elasticity with an imposed body load that is quadratic in $x_0$, the distance along the bar. The problem setup and boundary conditions are illustrated in Fig. \ref{Ex2_fig1} with $L = 1.0\, \mbox{m}$, $l = 0.5\, \mbox{m}$, and a cross-section area $A = l^2 = 0.25\, \mbox{m}^2$. The bar is made of a single linear elastic material with a Young's modulus $E = 10.0\, \mbox{N/m}^2$ and a Poisson ratio set to $\nu = 0.0$ to avoid three-dimensional effects. The bar is clamped at its left extremity, and $\mathbf{u}_D = [u_{Dx_0}\ u_{Dy_0}\ u_{Dz_0}]^T = [0.0\ 0.0 \ 0.0]\, \mbox{m}$. For all simulations, the Nitsche's penalty parameter is set to $\gamma_N = 100.0$. When applying the ghost stabilization, the penalty parameter is fixed to $\gamma_G = 0.001$. A quadratic body load $b = [b_{x_0}\ b_{y_0}\ b_{z_0}]^T = b_0\, [x_0^2\ 0.0\ 0.0]\, \mbox{N/m}^3$ with $b_{x_0} = 2.0$ is applied along the bar. This loading case yields a quartic displacement solution over the bar: 
\begin{equation}
\begin{array}{lll}
u_{x_0}(x) & = & \displaystyle u_{Dx_0} + \frac{b_0}{12EA} \left( 4L^3x_0 - x_0^4\right),\ \mbox{with}\quad b_{x_0} = b_0 x_0^2,\\[2.5pt]
u_{y_0} & = & 0,\\[5pt]
u_{z_0} & = & 0.
\end{array}
\label{Ex2_Eq1_analyticalSolution}
\end{equation}

\begin{figure*}[!h]\center
\includegraphics[scale=1]{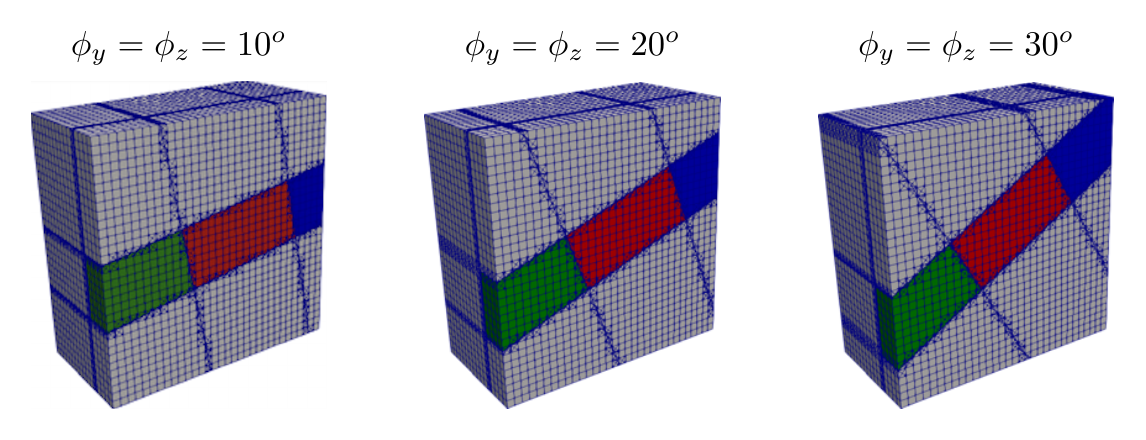}
\caption{Immersed bar rotated with an angle $\phi_y$ around the y-axis and an angle $\phi_z$ around the z-axis in a fixed three dimensional background mesh.}
\label{Ex2_fig2}
\end{figure*}

The bar is rotated by an angle $\phi_y$ around the y-axis and $\phi_z$ around the z-axis. The following angles with respect to the orientation of the background mesh are considered $\phi_y = \phi_z = [\, 10^o\ 20^o\ 30^o\ 40^o\ 50^o\ 60^o\ 70^o\ 80^o\, ]$, as illustrated in Fig. \ref{Ex2_fig2}. The background mesh size is successively refined, and the following mesh sizes are considered: $h = [0.5\ 0.25\ 0.125\ 0.0625]\, \mbox{m}$. Each setup is solved with linear, quadratic, and cubic B-splines. The system condition number, as defined in Eq.~\eqref{Eq_Cond}, is monitored. The accuracy of the evaluated physical responses is compared against the analytical solution given in Eq.~\eqref{Ex2_Eq1_analyticalSolution} using the relative $L^2$ error norm and the relative $H^1$ error semi-norm, as defined in Eqs.~(\ref{Eq_L2}, \ref{Eq_H1}).

\begin{figure*}[!h]\center
\includegraphics[scale=0.75]{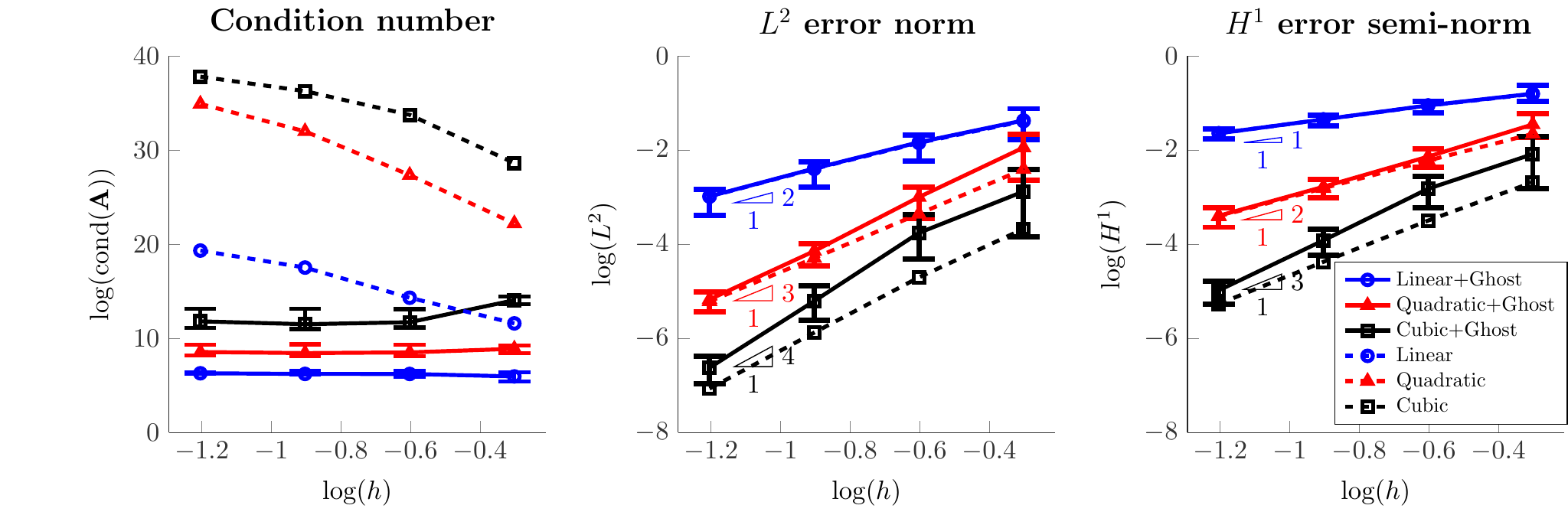}
\caption{Robustness study on an immersed elastic bar in a fixed 3D background mesh: mean performance and error over eight rotation angles.}
\label{Ex2_fig3}
\end{figure*}

For the configurations defined above, the condition number, the $L^2$ error norm, and the $H^1$ error semi-norm are given in Fig. \ref{Ex2_fig3}. The first graph shows the evolution of the mean condition number averaged over the eight rotation angles in terms of the background mesh size. The second and third graphs present the mean relative $L^2$ error norm and the mean relative $H^1$ error semi-norm over the eight rotation angles with respect to the background mesh size. For each performance measure, the solid and dashed lines correspond to the application or the absence of the ghost penalty stabilization respectively. The error bars illustrate the range of $L^2$ and $H^1$ values when using the ghost penalty. The range is defined by the minimum and maximum values of the relative $L^2$ and $H^1$ errors. The use of linear, quadratic, or cubic B-splines is indicated by circle, triangle, or square marks.

When no stabilization is used (dashed lines in Fig. \ref{Ex2_fig3}), an increase of the condition number is observed when using higher-order basis functions, and when using finer meshes. Introducing ghost penalty stabilization (solid lines in Fig. \ref{Ex2_fig3}), the condition number can be significantly improved, especially when using higher-order basis functions, and the effect of the mesh refinement is mitigated. A higher accuracy is achieved when using higher order B-splines, as expected since the exact solution is quartic along the bar. The relative $L^2$ error norm and the relative $H^1$ error semi-norm values are similar with and without ghost stabilization. It should be noted that for cubic B-splines, ghost stabilization leads to a drastic improvement of the condition number at the price of a slight increase of the $L^2$ and $H^1$ errors. The error bars show that, using ghost stabilization, the finite element predictions are insensitive to the rotation of the bar. This suggests that the proposed XIGA approach exhibits robustness with respect to intersection configurations. Finally, the graphs show that the optimal convergence rates with respect to mesh refinement are recovered for all B-spline orders $p$, i.e., $p+1$ in the $L^2$ error norm and $p$ in the $H^1$ error semi-norm.

\subsection{Accuracy study for \textit{N}-phase problems}\label{Ex3}
\begin{figure*}[!h]\center
\includegraphics[scale=1]{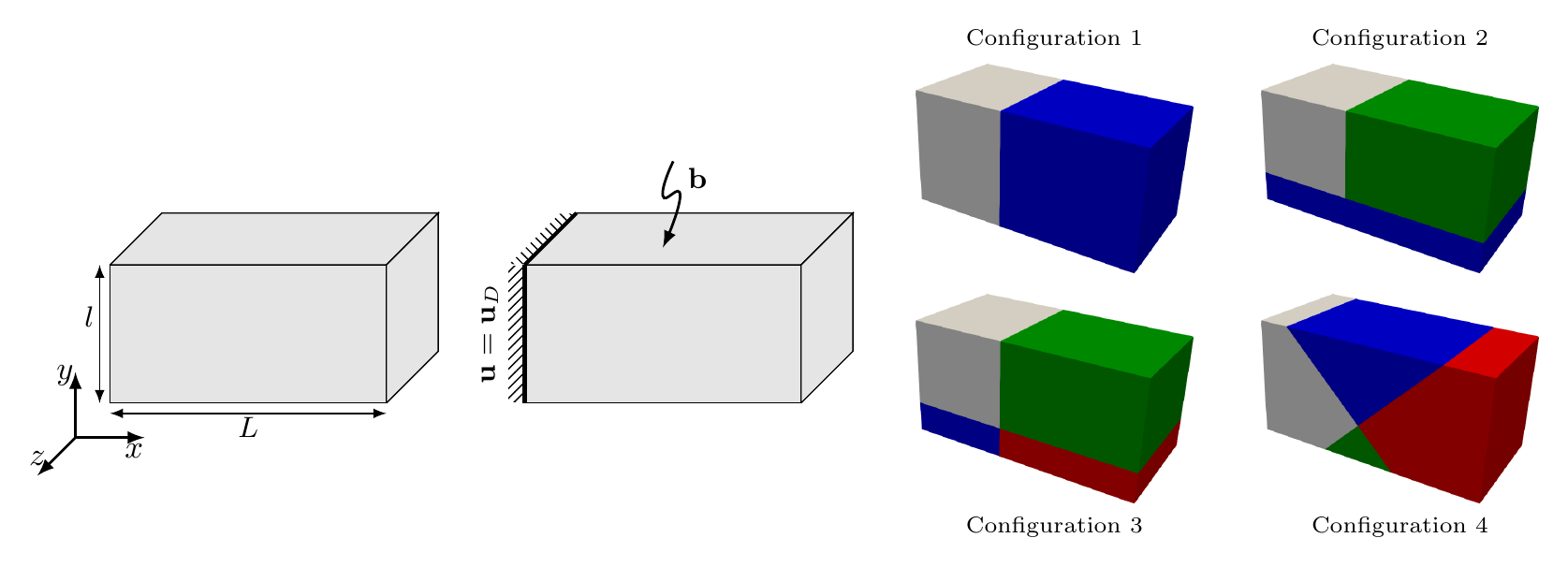}
\caption{Problem setup, boundary conditions, and configurations for a linear elastic single material bar with \textit{N}-phase junction.}
\label{Ex3_fig1}
\end{figure*}

In this section, the ability of the XIGA approach to handle \textit{N}-phase junctions is investigated. Again we consider a linear elastic problem with an imposed body load. The problem setup and boundary conditions are illustrated in Fig. \ref{Ex3_fig1} with $L = 1.0\, \mbox{m}$, $l = 0.5\, \mbox{m}$, and a cross-section area $A = l^2 = 0.25\, \mbox{m}^2$. The bar problem presented in Subsection \ref{Ex2} is reused, but with the geometric configurations illustrated in Fig. \ref{Ex3_fig1}. The bar is made of four different phases filled with the same material, so that the numerical solution can be easily compared to an analytical one. The Young's modulus of the material is set to $E = 10.0\, \mbox{N/m}^2$ and a Poisson ratio $\nu = 0.0$ is chosen to avoid three-dimensional effects. The bar is clamped at its left extremity, $\mathbf{u}_D = [u_{Dx}\ u_{Dy}\ u_{Dz}]^T = [0.0\ 0.0 \ 0.0]\, \mbox{m}$.  For all simulations, the Nitsche's penalty parameter is set to $\gamma_N = 100.0$ and the ghost penalty parameter is set to $\gamma_G = 0.001$. A quadratic body load $b = b_0\, [x^2\ 0.0\ 0.0]\, \mbox{N/m}^3$ with $b_0 = 2.0$ is imposed, leading to a quartic displacement solution over the bar, as described hereunder:
\begin{equation}
\begin{array}{lll}
u_x(x) & = & \displaystyle u_{Dx} + \frac{b_0}{12EA} \left( 4L^3x - x^4\right),\ \mbox{with}\quad b_x = b_0 x^2,\\[2.5pt]
u_y & = & 0,\\[5pt]
u_z & = & 0.
\end{array}
\label{Ex3_Eq1_analyticalSolution}
\end{equation}

As illustrated in Fig. \ref{Ex3_fig1}, four different configurations are considered: (1) a two-phase junction with the interface aligned with the $y$-axis, (2) a three-phase junction with the interfaces aligned with the $x$- and $y$-axes, (3) a four-phase junction with the interfaces aligned with the $x$- and $y$-axes, and (4) a four-phase junction with the interfaces rotated by 45 degrees against with the $x$- and $y$-axes. An \textit{h}-refinement study is performed to investigate the convergence rates of linear, quadratic, and cubic B-spline basis functions. Five mesh sizes are considered, $h = [ 0.5\ 0.25\ 0.125\ 0.0625\ 0.03125]\, \mbox{m}$.

\begin{figure*}[!h]\center
\includegraphics[scale=0.75]{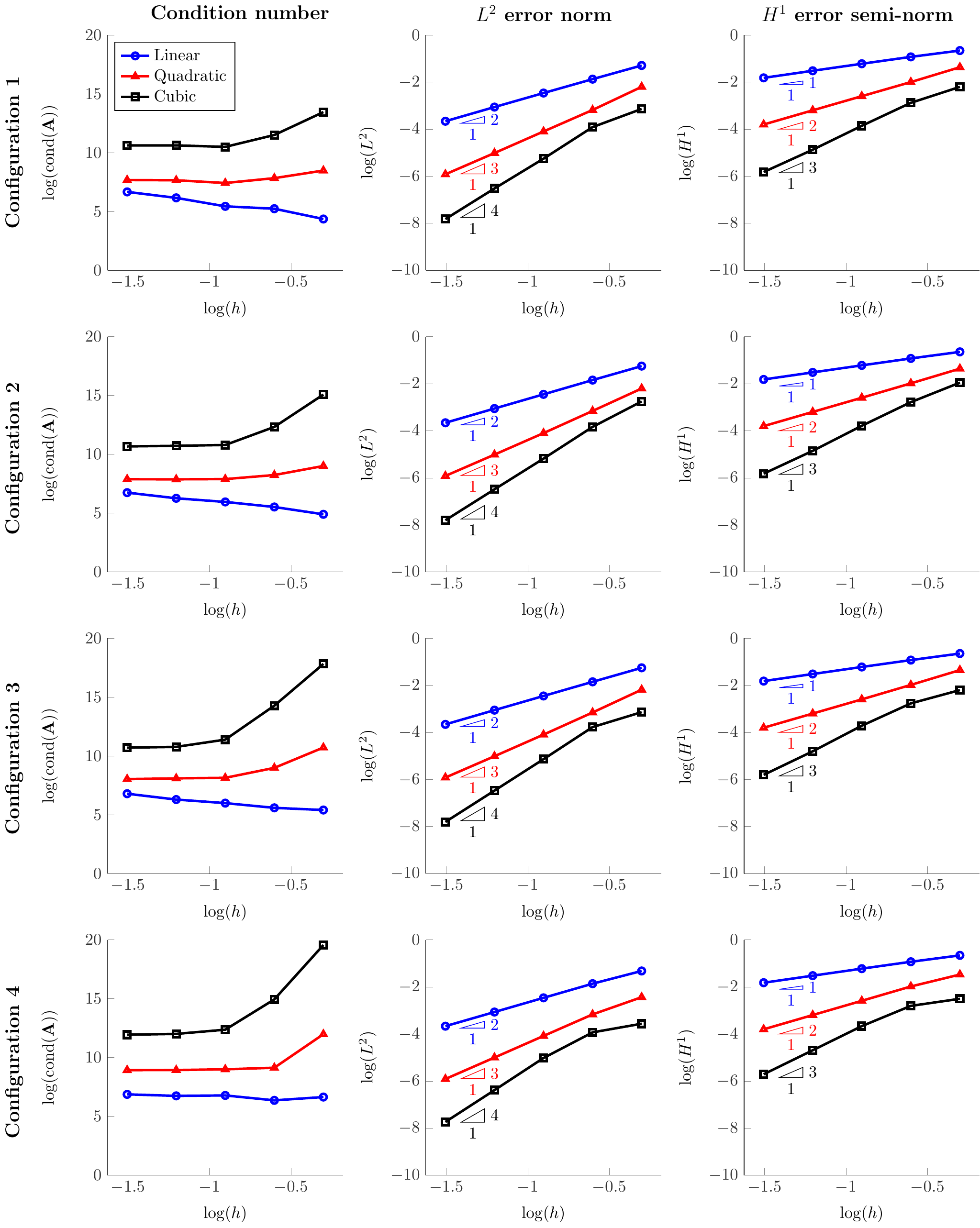}
\caption{Accuracy study on a linear elastic single material bar with \textit{N}-phase junction configurations in three dimensions.}
\label{Ex3_fig2}
\end{figure*}

The problem is solved for the four geometric configurations depicted in Fig. \ref{Ex3_fig1}, successively refining the mesh. The performance in terms of system conditioning and solution accuracy is provided in Fig. \ref{Ex3_fig2}. The graphs show the variations of the condition number, the relative $L^2$ error norm and the relative $H^1$ error semi-norm for different mesh sizes. For each performance measure, the use of linear, quadratic, or cubic B-splines is indicated by circle, triangle, or square marks.

For all configurations, similar performance is observed in terms of condition number, $L^2$, and $H^1$ errors. This demonstrates the ability of the approach to handle \textit{N}-phase junctions accurately. As expected, an increase in the condition number is observed for higher-order B-splines. In terms of accuracy, optimal convergence rates are recovered for both the $L^2$ and the $H^1$ error norms for all basis function orders. It should be noted that a slight increase in the condition number and in the errors is observed as the complexity of the geometric configurations is increased. Additionally, an increase in the condition number is observed for low refinement, due to the creation of small material integration subdomains with low volume ratio with respect to the background elements. 

\subsection{Accuracy study for two-material problems with curved interface}\label{Ex4}
\begin{figure}[h!]\center
\includegraphics[scale=1]{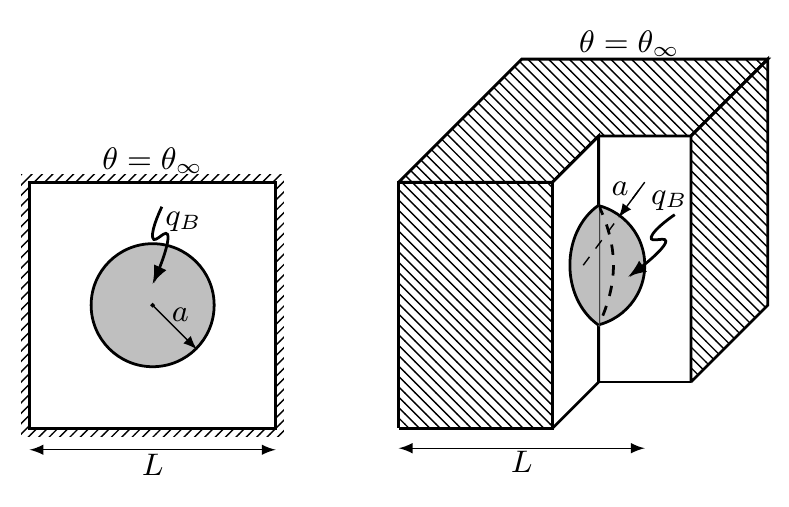}
\caption{Problem setup and boundary conditions for a heated inclusion in an infinite host medium in two and three dimensions.}
\label{Ex4_fig1}
\end{figure}

In the previous subsection, the surfaces and interfaces were planar and could be represented exactly with our level set approach, that inherently leads to a low order approximation of geometry, see Subsection \ref{sectionGeometry}. In this subsection, the ability of the proposed XIGA approach to accurately resolve non-planar surfaces and interfaces is studied. A heated inclusion embedded in an infinite matrix problem is considered in two and three dimensions. The problem setup is illustrated in Fig. \ref{Ex4_fig1}, where the dimensions are set to $L = 2.0\, \mbox{m}$ and $a = 0.5\, \mbox{m}$. The embedded inclusion is made of a material $I$ with a conductivity $\kappa^I = 1.0\, \mbox{W/mK}$ and is undergoing a constant heat body load $q_B = 1.0\, \mbox{W/m}^2$ or $\mbox{W/m}^3$. The infinite medium is made of a material $II$ with a conductivity $\kappa^{II} = 0.125\, \mbox{W/mK}$. The surrounding medium is not heated, i.e, $q_B^{II} = 0.0\, \mbox{W/m}^2$ or $\mbox{W/m}^3$. For the 2D case, a cylinder embedded in an infinite medium is considered. In cylindrical coordinates, the strong form of the heat conduction equation is given as:
\begin{equation}
\frac{1}{r} \frac{d}{dr} \left( r \frac{d \theta(r)}{dr} \right) + \frac{q_B}{\kappa} = 0,
\label{Ex4_Eq1_diffEquation}
\end{equation}
where $r$ is the radius computed from the center of the cylinder and has the general solution:
\begin{equation}
\theta(r) = -\frac{q_B r^2}{4 \kappa} + C_1 \ln(r) + C_2,
\label{Ex4_Eq1_genSolEquation}
\end{equation}
where $C_1$ and $C_2$ are integration constants that can be determined by considering the boundary conditions. A prescribed temperature $\theta_D = 0.375\, \mbox{K}$ and a zero temperature gradient are prescribed at the origin. At the material interface, continuity is enforced for the temperature field and the heat flux. Finally, the temperature field over the computational domain is given as:
\begin{equation}
\theta(r) = \left\{
\begin{array}{lll}
\displaystyle \theta_D -\frac{q_B^I r^2}{4 \kappa^I},                                                    &\mbox{for}& r \leq a,\\[7.5pt]
\displaystyle \theta_D -\frac{q_B^I a^2}{4 \kappa^I} - \frac{q_B^I a^2}{2 \kappa^{II}} \ln\left(\frac{r}{a}\right), &\mbox{for}& r > a.
\end{array}
\right.
\label{Ex4_Eq1_diffEquation}
\end{equation}
For the 3D heat conduction problem, a sphere embedded in an infinite medium is considered. In spherical coordinates, the strong form of the heat conduction equation is given as:
\begin{equation}
\frac{1}{r^2}\frac{d}{dr} \left( r^2 \frac{d \theta(r)}{dr} \right) + \frac{q_B}{\kappa} = 0,
\label{Ex4_Eq2_diffEquation}
\end{equation}
where $r$ is the radial coordinate. The general solution of Eq.~\eqref{Ex4_Eq2_diffEquation} is:
\begin{equation}
\theta(r) = -\frac{q_B r^2}{6 \kappa} + \frac{C_1}{r} + C_2,
\label{Ex4_Eq2_genSolEquation}
\end{equation}
where $C_1$ and $C_2$ are integration constants that can be determined by considering the boundary conditions. A prescribed temperature $\theta_D = 0.375\, \mbox{K}$ and a zero temperature gradient are prescribed at the origin. At the material interface, continuity is enforced for the temperature field and the heat flux. Finally, the temperature field over the computational domain is given as:
\begin{equation}
\theta(r) = \left\{
\begin{array}{lll}
\displaystyle \theta_D - \frac{q_B^I r^2}{6 \kappa^I},                                                     &\mbox{for}& r \leq a,\\[7.5pt]
\displaystyle \theta_D - \frac{q_B^I a^2}{6 \kappa^I} - \frac{q_B^I a^3}{3 \kappa^{II}} \Big( \frac{1}{a} - \frac{1}{r} \Big), &\mbox{for}& r > a.
\end{array}
\right.
\label{Ex4_Eq2_diffEquation}
\end{equation}
To simulate an infinite host domain, the analytical solutions, presented in Eqs.~(\ref{Ex4_Eq1_diffEquation}, \ref{Ex4_Eq2_diffEquation}), are applied to the outer faces of the host domain as a weak Dirichlet boundary condition. The Nitsche's penalty parameter is set to $\gamma_N = 100.0$ for both the boundary and the interface conditions; the ghost penalty parameter is set to $\gamma_G = 0.001$. Simulations are performed for different mesh sizes $h = [0.5\ 0.25\ 0.125\ 0.06125\ 0.03125]\, \mbox{m}$. The accuracy and convergence of the physical responses with mesh refinement are evaluated using the relative $L^2$ error norm and the relative $H^1$ error semi-norm, as defined in Eqs.~(\ref{Eq_L2}, \ref{Eq_H1}). 

\begin{figure*}[!h]\center
\subfigure[][$h=0.5\, \mbox{m}$.]{\includegraphics[width=2.75cm]{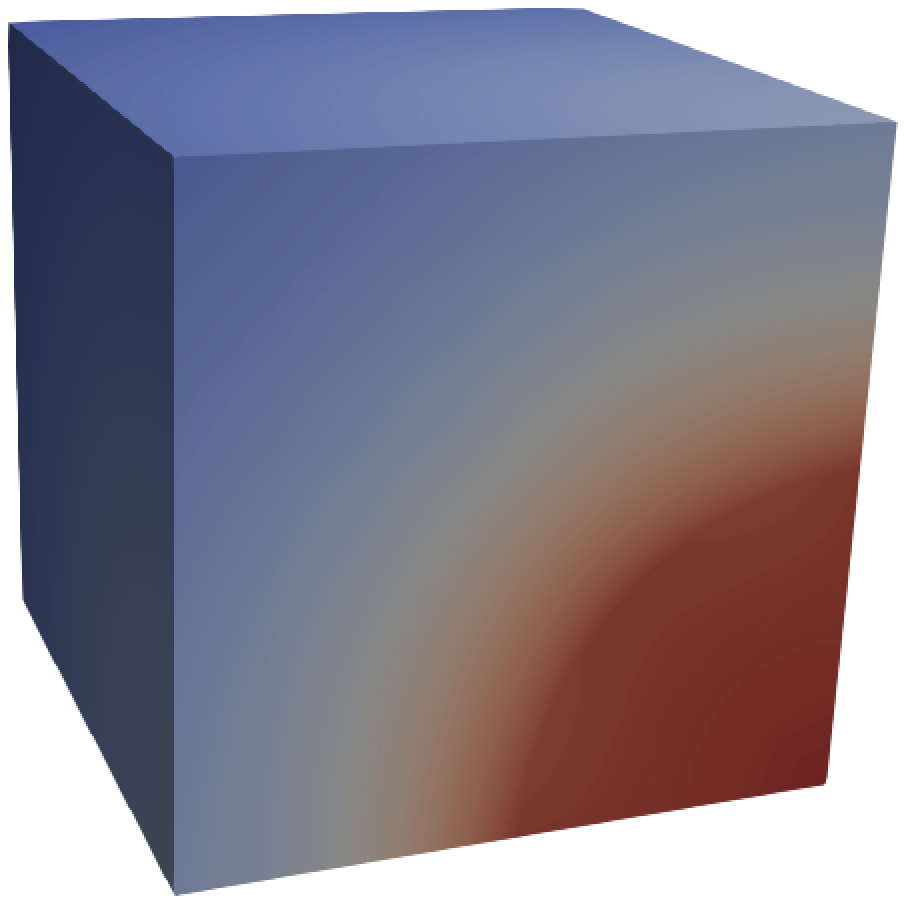}}\hskip0.2cm
\subfigure[][$h=0.25\, \mbox{m}$.]{\includegraphics[width=2.75cm]{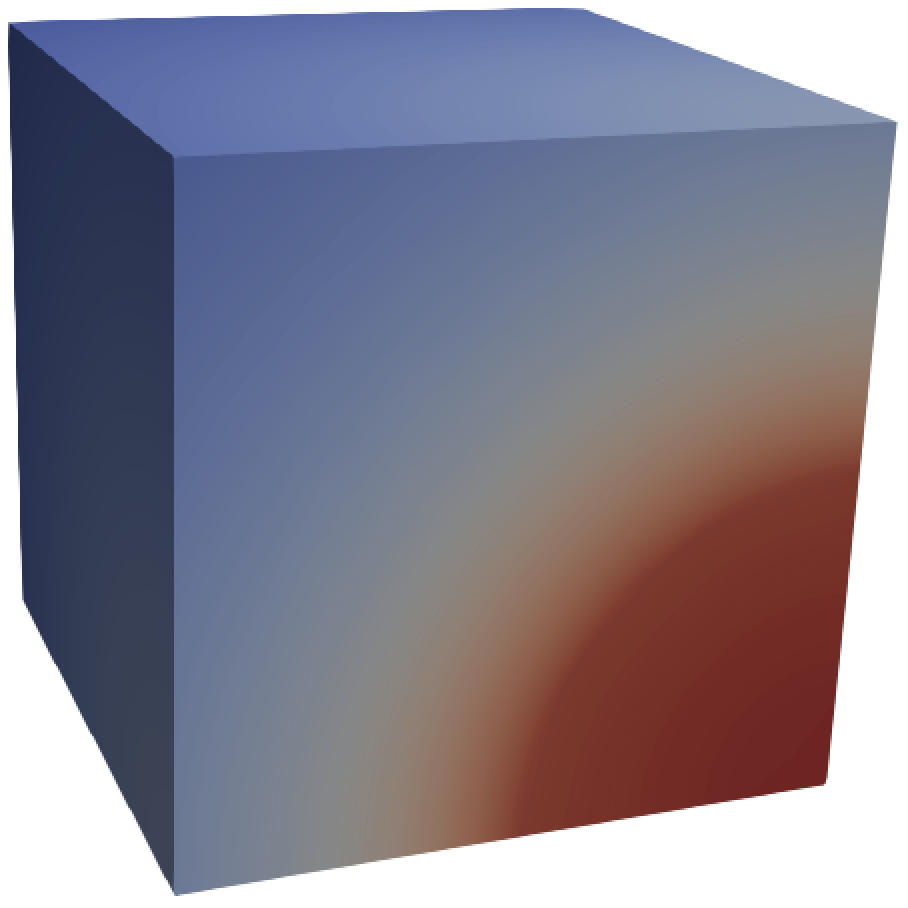}}\hskip0.2cm
\subfigure[][$h=0.125\, \mbox{m}$.]{\includegraphics[width=2.75cm]{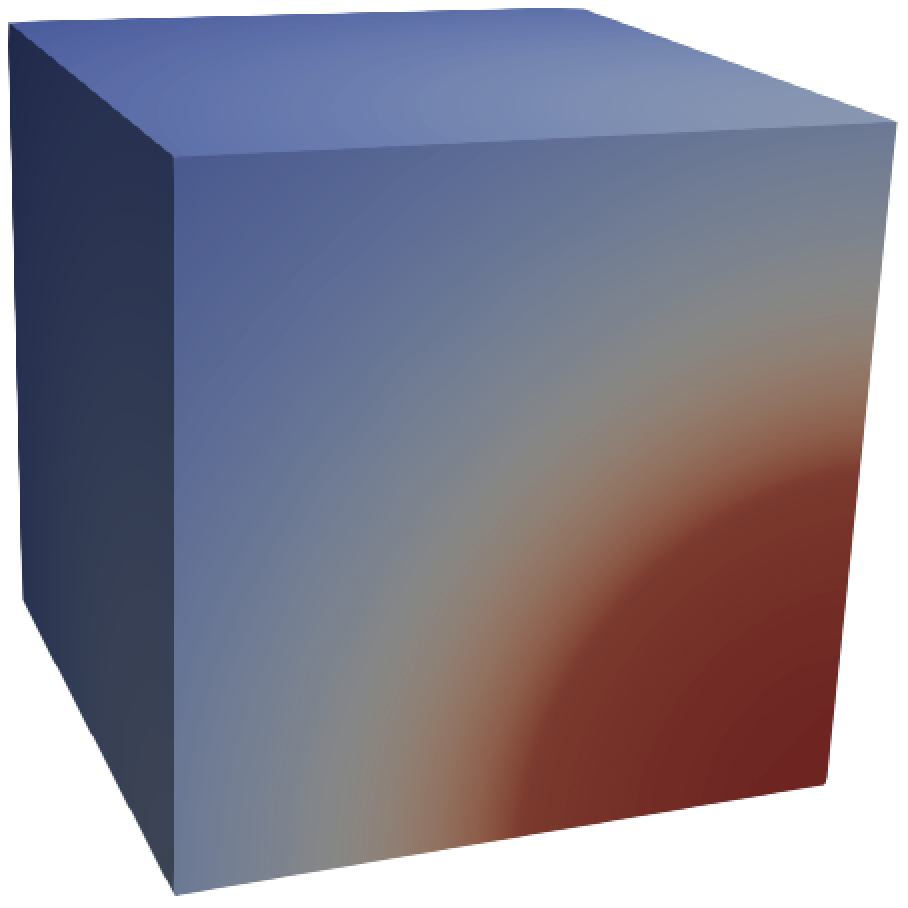}}\hskip0.2cm
\subfigure[][$h=0.0625\, \mbox{m}$.]{\includegraphics[width=2.75cm]{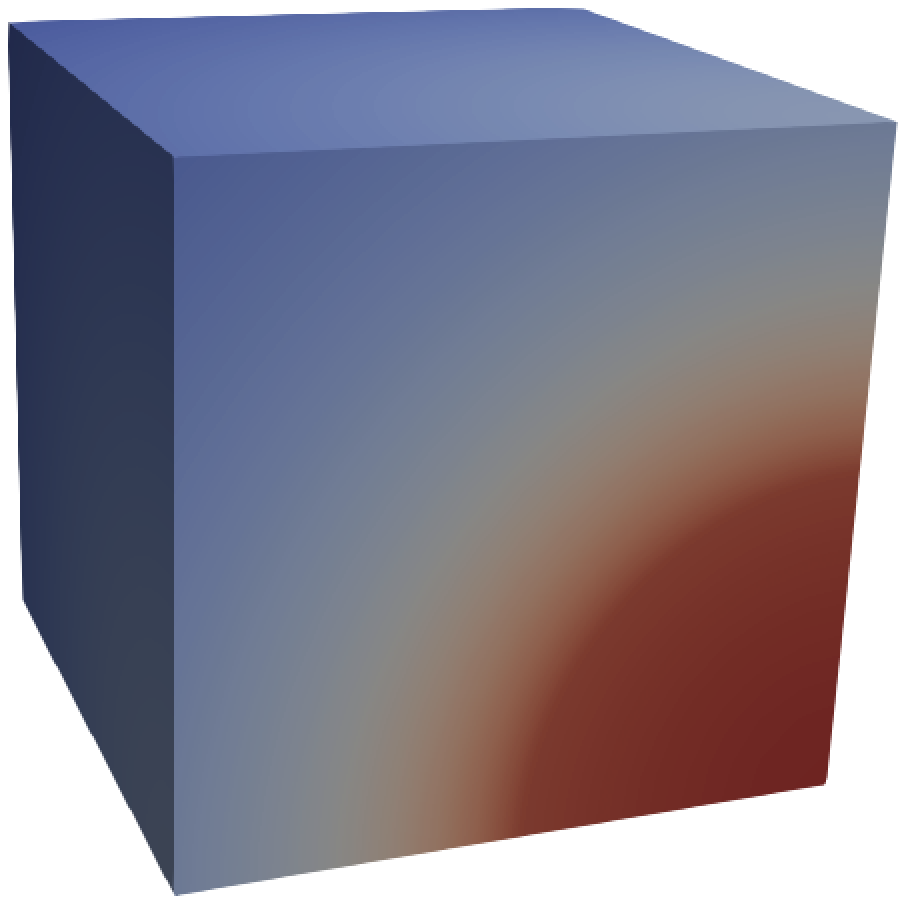}}\hskip0.2cm
\subfigure[][$h=0.03125\, \mbox{m}$.]{\includegraphics[width=2.75cm]{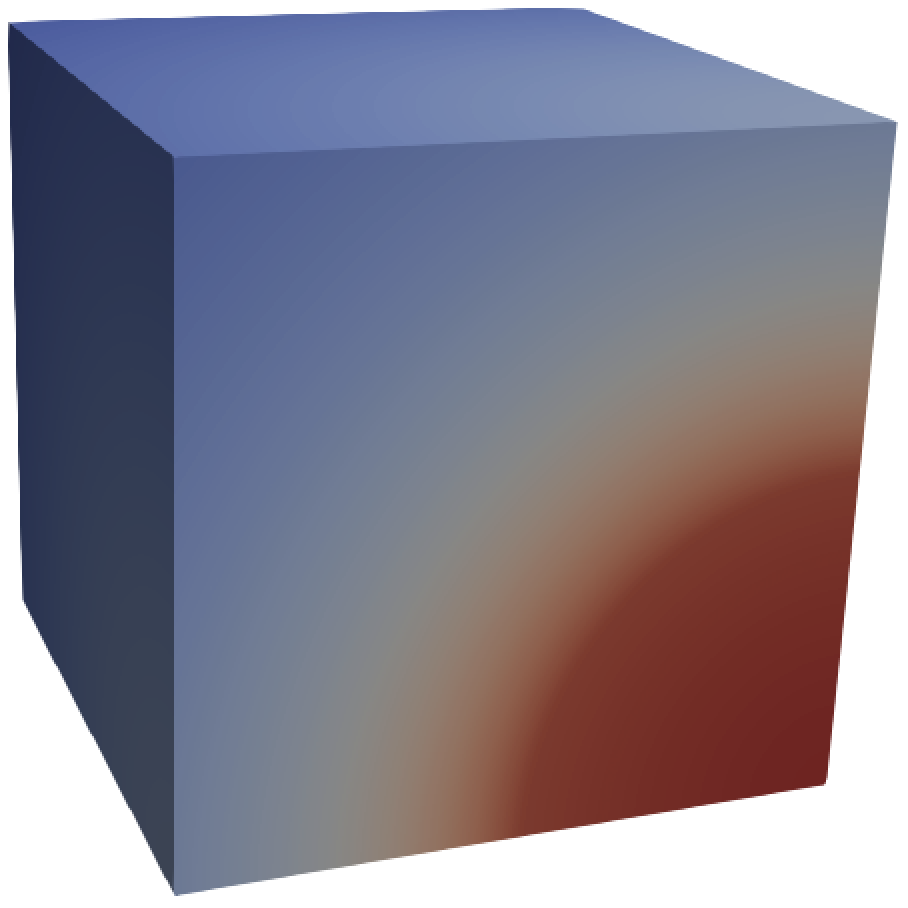}}\hskip0.2cm
\includegraphics[scale=1.0]{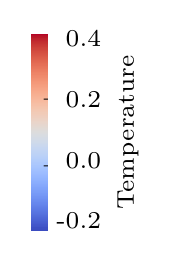}
\caption{Heated inclusion in an infinite medium, refinement of the temperature field with a fixed geometry field $h_{\mbox{\tiny int}} = 0.03125\, \mbox{m}$, solution using quadratic B-splines in three dimensions.}
\label{Ex4_fig2}
\end{figure*}

The temperature solution obtained for the five different mesh sizes with quadratic B-splines and a fixed integration mesh of size $h_{\mbox{\tiny int}} = 0.03125\, \mbox{m}$ is shown in Fig. \ref{Ex4_fig2}. Note that the geometry and the temperature fields are refined independently. The geometry is linearly interpolated and its accuracy can be increased by refining the integration mesh. A geometrical error, characterizing the accuracy of the cylinder or sphere representation, is evaluated as follows:
\begin{equation}
e_{\mbox{\tiny geo}} = \frac{V^h - V}{V},
\label{Ex4_Eq3_geoError}
\end{equation}
where $V$ is the reference volume or area in 2D or 3D, here chosen as $V = \pi a^2$ and $V = 4 \pi a^3/3$ in 2D and 3D respectively, and $V^h$, the numerical volume or area evaluated with the XIGA approach.

\begin{figure*}[!h]\center
\includegraphics[scale=0.75]{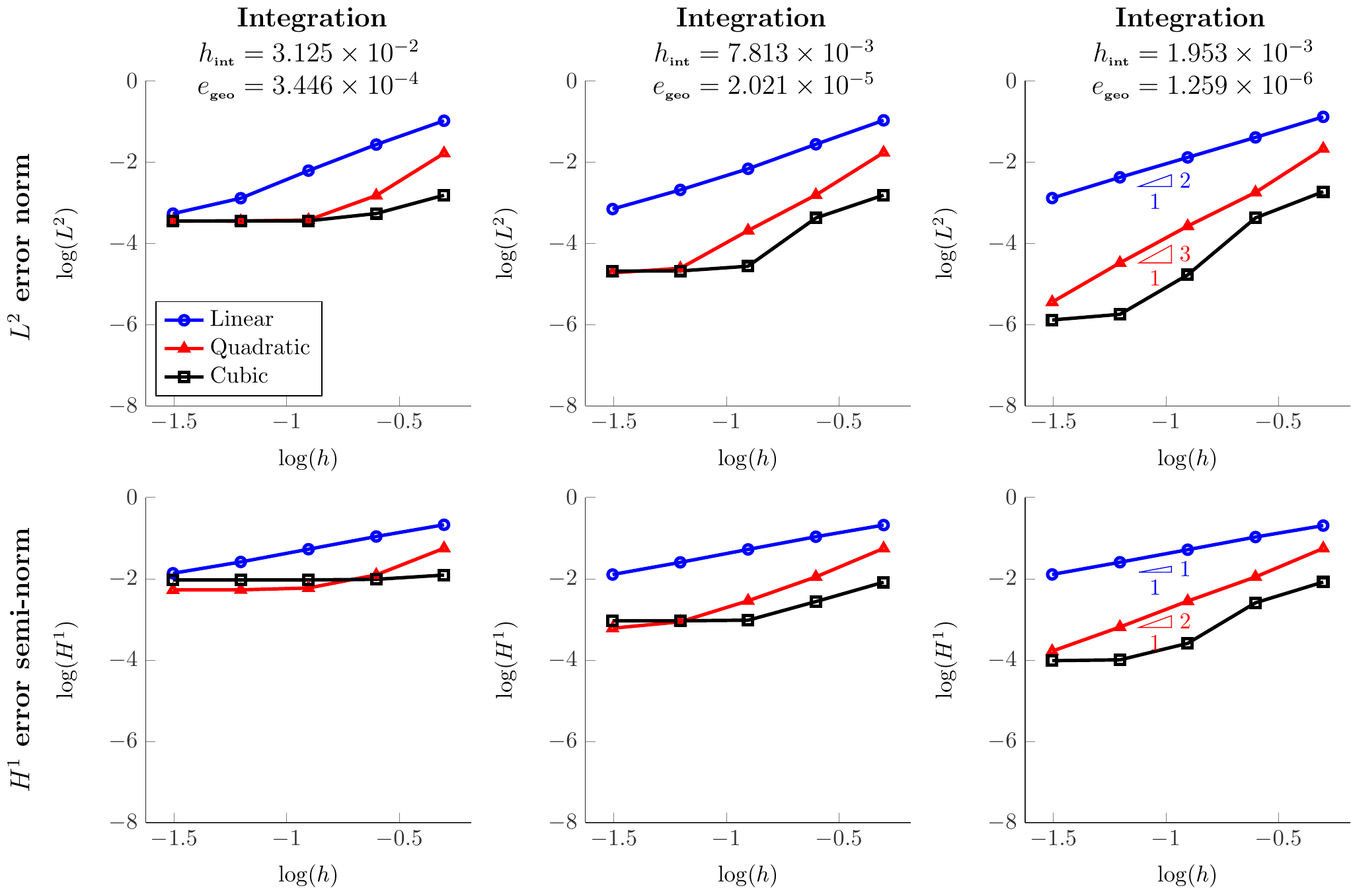}
\caption{Heated cylindrical inclusion in an infinite medium, refinement of the state field for three different fixed geometry refinements $h_{\mbox{\tiny int}} = 3.125e^{-2}$, $7.813e^{-3}$, and $1.953e^{-3}\, \mbox{m}$.}
\label{Ex4_fig3a}
\end{figure*}

\begin{figure*}[!h]\center
\includegraphics[scale=0.75]{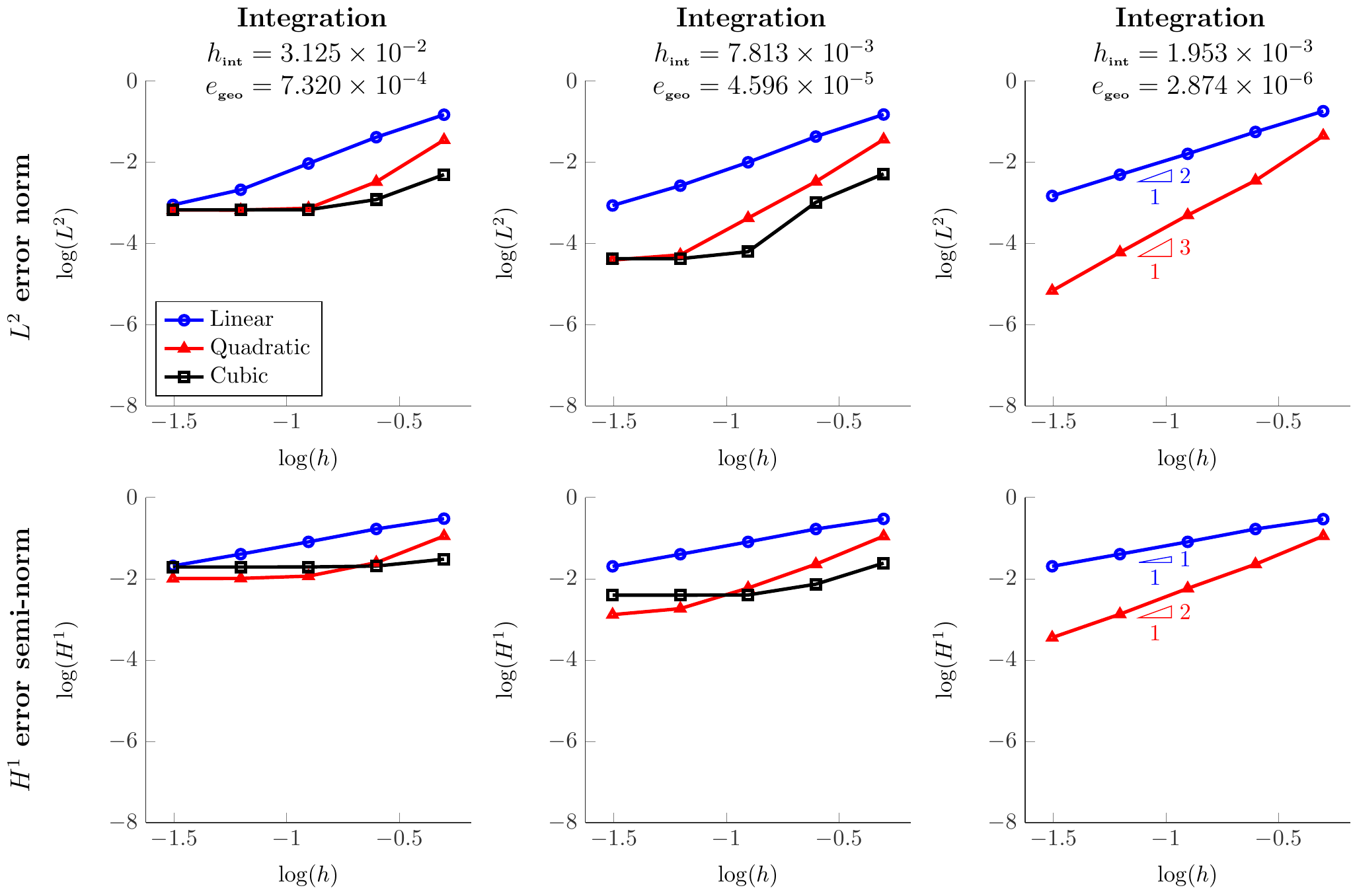}
\caption{Heated spherical inclusion in an infinite medium, refinement of the state field for three different fixed geometry refinements $h_{\mbox{\tiny int}} = 3.125e^{-2}$, $7.813e^{-3}$, and $1.953e^{-3}\, \mbox{m}$.}
\label{Ex4_fig3b}
\end{figure*}

The results obtained in two and three dimensions, are shown in Fig. \ref{Ex4_fig3a} and \ref{Ex4_fig3b}. In these figures, the rows present the relative $L^2$ error norm and the relative $H^1$ error semi-norm with respect to the \textit{h}-refinement. The columns correspond to different levels of refinement of the integration mesh, i.e., $h_{\mbox{\tiny int}} = 3.125e^{-2}$, $7.813e^{-3}$, and $1.953e^{-3}\, \mbox{m}$, that yield decreasing geometrical error $e_{\mbox{\tiny geo}}$, as shown on top of each graph. The highest level of refinement leads to a problem with around 12,200 integration elements and around 1,300 DOFs in 2D and around 11,500,000 integration elements and around 41,000 DOFs in 3D. Therefore, the cubic case with the highest refinement of the integration mesh is omitted in 3D, as the computational resources to tackle such a problem were not available.

In 2D and in 3D, the relative $L^2$ and $H^1$ error plots show that the ability of the proposed approach to represent the curved two-material interface problem solutions increases as the B-spline mesh is refined. However, optimal convergence rates in $L^2$ and $H^1$ errors are not fully recovered for all B-splines orders and geometry representations. When the integration mesh is not sufficiently refined, using higher-order bases only provides a slight improvement over using a linear basis. In this case, the lack of a sufficient geometry resolution limits the accuracy of the computed solution and prevents the recovery of optimal convergence rates of the $L^2$ and $H^1$ errors with \textit{h}-refinement. As the integration mesh is further refined, the geometry error on the curved interface of the circle or the sphere drops and this issue is mitigated. Using an integration mesh size of $h_{\mbox{\tiny int}} =1.953e^{-3}\, \mbox{m}$, optimal convergence rates in the relative $L^2$ and $H^1$ error norm is fully recovered for the quadratic B-splines and partially recovered for the cubic ones. 

\subsection{Accuracy study for \textit{N}-material problems}\label{Ex5}
\begin{figure}[h!]\center
\includegraphics[scale=1]{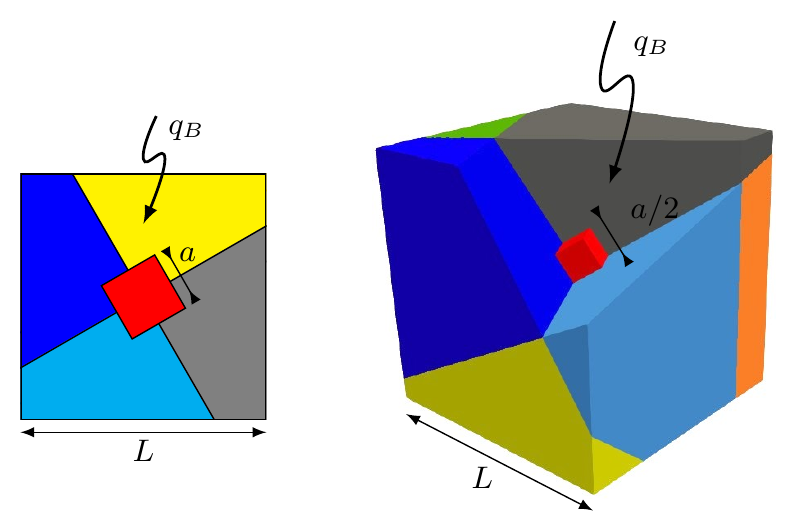}
\caption{Setup and boundary conditions for a heated multi-material medium in two and three dimensions.}
\label{Ex5_fig1}
\end{figure}

Finally, the ability to accurately tackle multi-material problems in two and three dimensions is addressed in this subsection. In two dimensions, the case of a beam with a square cross-section embedded in a four-material host matrix is considered. In three dimensions, a cubic inclusion is embedded in an eight-material matrix. The setup of the problem is illustrated in Fig. \ref{Ex5_fig1}, where the dimensions are set to $L = 2.0\, \mbox{m}$ and $a = 0.5\, \mbox{m}$. A fixed temperature $\theta = 0\, \, \mbox{K}$ is imposed on the outer faces of the host domain as a weak Dirichlet boundary condition. The ghost penalty parameter is set to $\gamma_G = 0.001$ and the Nitsche's penalty parameter to $\gamma_N = 50.0$ for both the boundary and the interface conditions. The entire domain, inclusion and matrix, is undergoing a varying heat body load $q_B = \sin(2\pi x) \sin(2\pi y)\, \mbox{W/m}^2$ in 2D and $\sin(2\pi x) \sin(2\pi y) \sin(2\pi z)\, \mbox{W/m}^3$ in 3D. Two different material configurations are considered: a single material one and a multi-material one. In the multi-material setting, singularities arise at the sharp corners of the inclusion. As this is not the case with a single material problem, the two configurations are considered to assess the effect of such singularities on the accuracy of the numerical solution. In the single material case, the inclusion and the four, in 2D, or eight, in 3D, matrix materials share the same conductivity $\kappa = 1.0\, \mbox{W/mK}$. In the multi-material case, the embedded inclusion is made of a material $I$ with a conductivity $\kappa^I = 1.0\, \mbox{W/mK}$ and the host medium is made of four, in 2D, or eight, in 3D, different materials with a conductivity $\kappa^{i} = i \times 0.125\, \mbox{W/mK}$.  

Simulations are performed using linear, quadratic, and cubic B-spline basis functions. The accuracy of the physical responses are evaluated by performing a \textit{h}-refinement with $h = [ 0.5\ 0.25\ 0.125\ 0.0625\ 0.03125]\, \mbox{m}$ and comparing the obtained results  in terms of the relative $L^2$ and $H^1$ error norms, as defined in Eqs.~(\ref{Eq_L2}, \ref{Eq_H1}), with the response evaluated with a higher resolution $h = 0.015625\,\mbox{m}$.

\begin{figure*}[!h]\center
\includegraphics[scale=1]{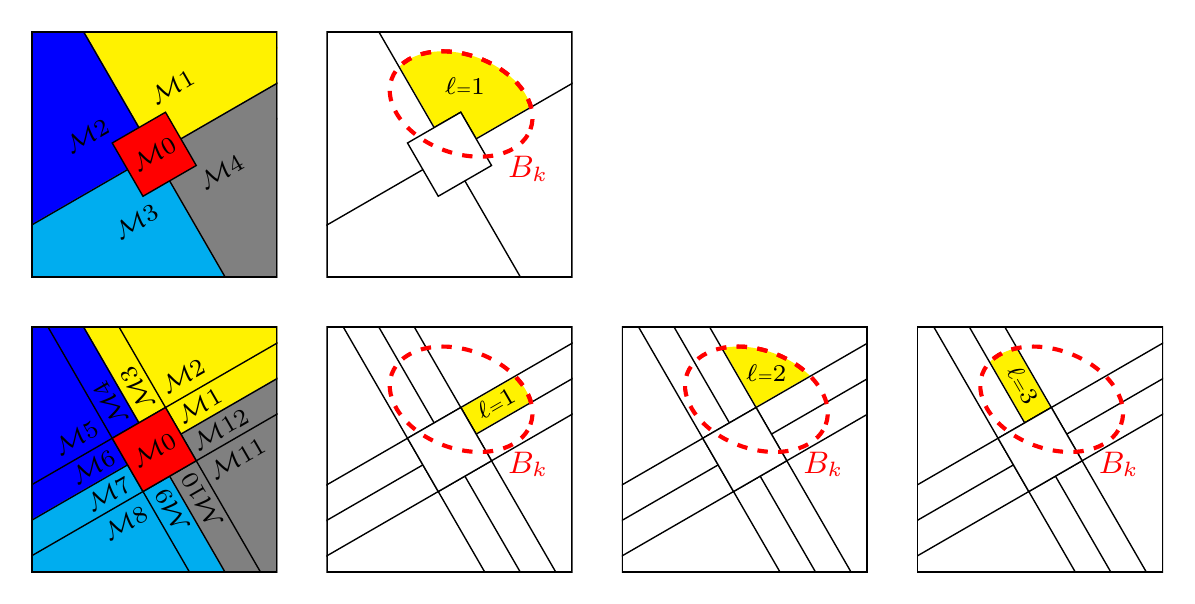}
\caption{Schemes for phase and material assignment for a heated single or multi-material medium.}
\label{Ex5_fig2}
\end{figure*}

The inclusion presents sharp corners which are with\-in a single background cell. Using the same basis function to approximate the temperature field around the corner might degrade the accuracy of the solution. To investigate this situation, two schemes for the phase and material assignments are considered for the two-dimensional case. Each region in the matrix associated with identical material properties can be divided into one or three connected subregions leading to a total number of five or thirteen materials, as illustrated in Fig. \ref{Ex5_fig2}. As detailed in Subsection \ref{subsectionEnrichment}, the two schemes lead to different enrichment. When considering five materials, a single enrichment level $\ell=1$ is necessary to approximate the temperature field around the inclusion corner; while for the thirteen-material case, three enrichment levels are used $\ell=1,\dots,3$. 

\begin{figure*}[!h]\center
\includegraphics[scale=0.75]{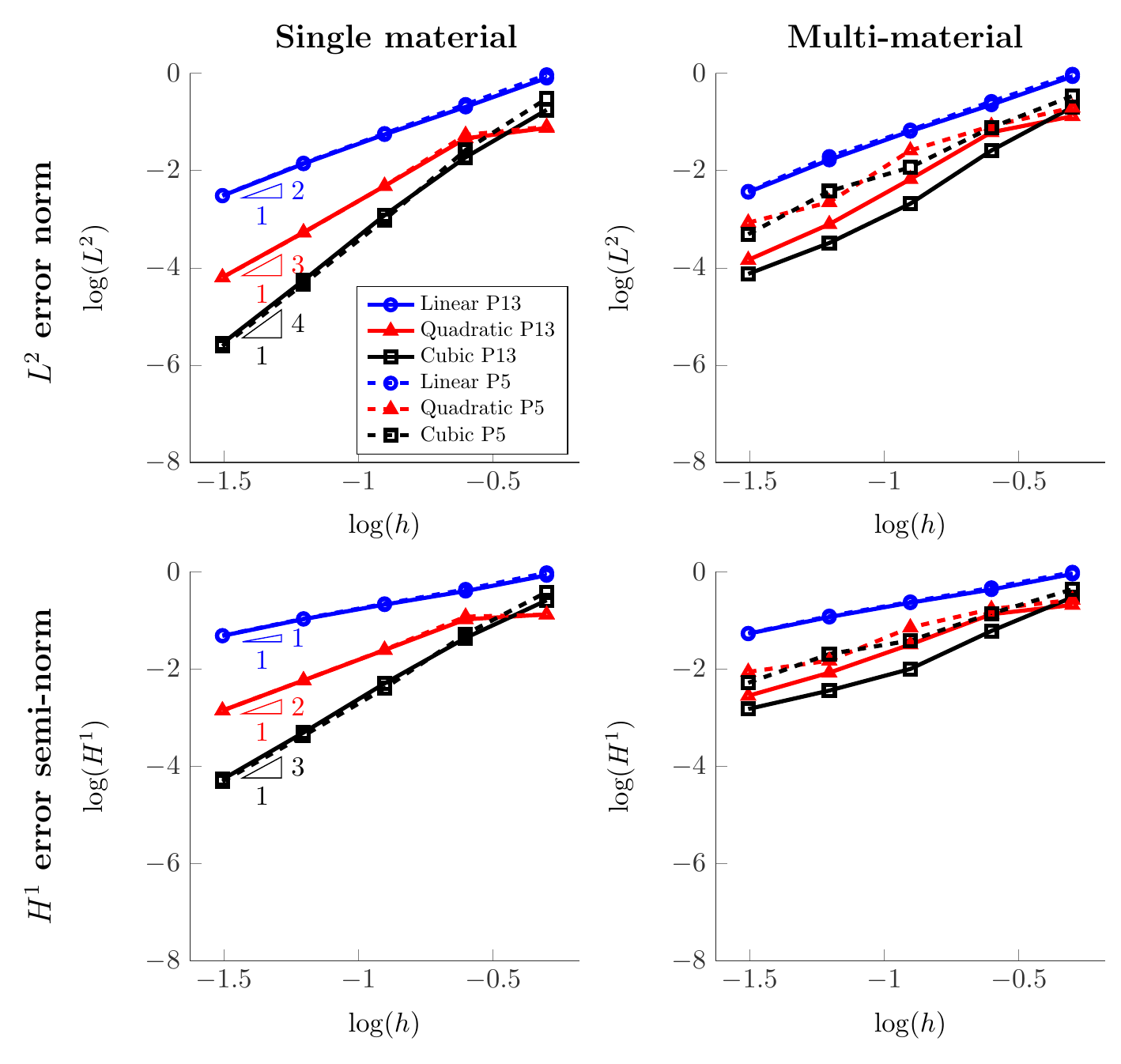}
\caption{Accuracy study on a heated single or multi-material medium in two dimensions.}
\label{Ex5_fig3}
\end{figure*}

For the two-dimensional case, the simulations use linear, quadratic, and cubic B-splines basis functions. The relative $L^2$ error norm and the relative $H^1$ error semi-norm with respect to the \textit{h}-refinement are presented in Fig. \ref{Ex5_fig3}. The highest level of refinement leads to a problem with 30,500 integration elements and about 22,000 DOFs. In the figure, the first column corresponds to the single material case and the second column to the multi-material case. Additionally, the results are presented for the five-material and the thirteen-material settings in solid and dashed line respectively. 

In 2D for the single material setting, optimal convergence rates in the relative $L^2$ and $H^1$ error norms are recovered for linear, quadratic, and cubic B-spline basis functions. For lower level of refinement, the solution obtained using cubic B-splines is less accurate than for quadratic B-splines. This lower accuracy results from the effective coarsening introduced by the ghost stabilization. For the multi-material setting, optimal convergence rates are only recovered for linear and quadratic B-splines basis functions. Using cubic basis functions provides a slight improvement in terms of accuracy, but not in terms of convergence rates with respect to the quadratic basis functions. This is due to singularities at the sharp corners of the inclusion. Introducing auxiliary material domains and enhancing the approximation around the inclusion sharp corner allows for the recovery of optimal convergence rates.    

\begin{figure*}[!h]\center
\subfigure[][$h=0.5\, \mbox{m}$.]{\includegraphics[width=2.75cm]{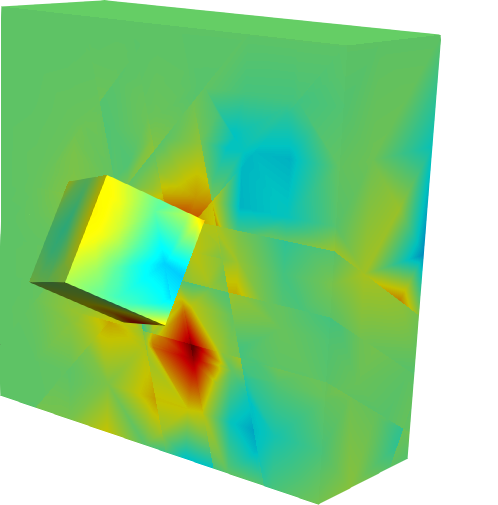}\label{Ex5_fig4a}}\hskip0.2cm
\subfigure[][$h=0.25\, \mbox{m}$.]{\includegraphics[width=2.75cm]{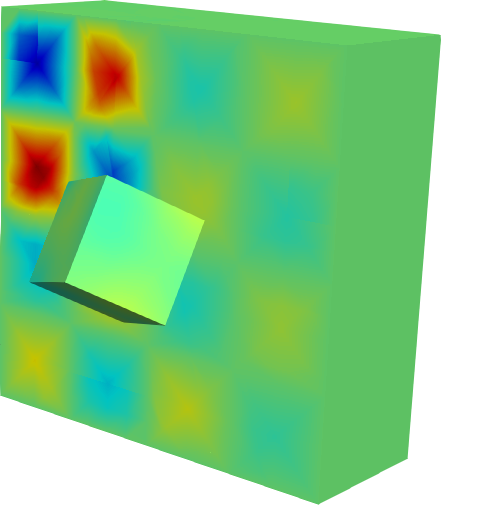}}\hskip0.2cm
\subfigure[][$h=0.125\, \mbox{m}$.]{\includegraphics[width=2.75cm]{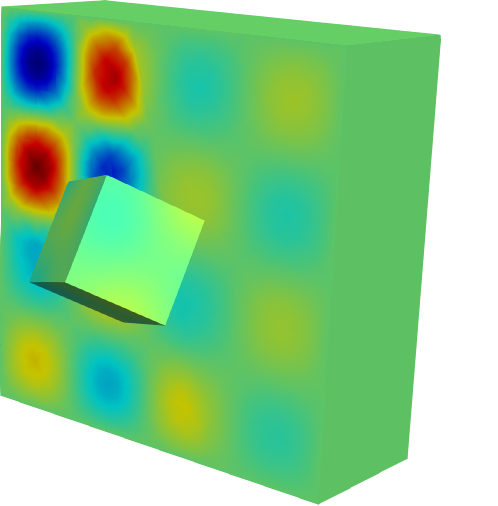}}\hskip0.2cm
\subfigure[][$h=0.0625\, \mbox{m}$.]{\includegraphics[width=2.75cm]{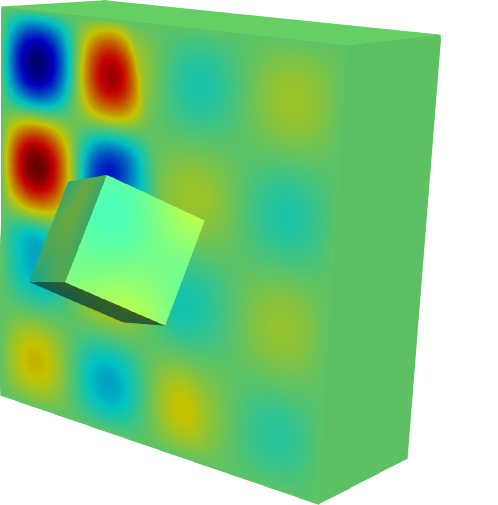}}\hskip0.2cm
\subfigure[][$h=0.03125\, \mbox{m}$.]{\includegraphics[width=2.75cm]{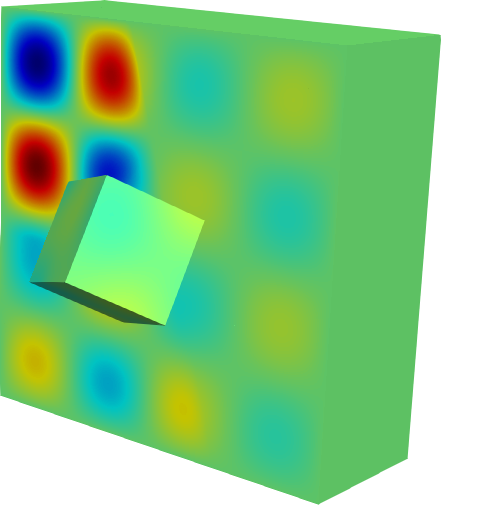}}\hskip-0.2cm
\includegraphics[scale=1]{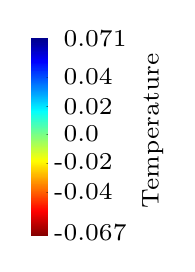}
\caption{Multi-material heated medium with refinement of the temperature field solution using quadratic B-splines in three dimensions.}
\label{Ex5_fig4}
\end{figure*}

For the three-dimensional case, only the scheme with the auxiliary material domains around the corners and edges is considered. The resulting setup has a total number of 57 materials; each area in the host medium is made of seven material domains with same properties. The simulations are carried out for quadratic B-splines basis functions only. The temperature solution obtained for the multi-material case and for the five different mesh sizes with quadratic B-spline is illustrated in Fig. \ref{Ex5_fig4}. It should be noted that for a coarse mesh size as shown in Fig. \ref{Ex5_fig4a} with $h=0.5\, \mbox{m}$, the nature of the solution is not representable by the spline space.

\begin{figure*}[!h]\center
\includegraphics[scale=0.75]{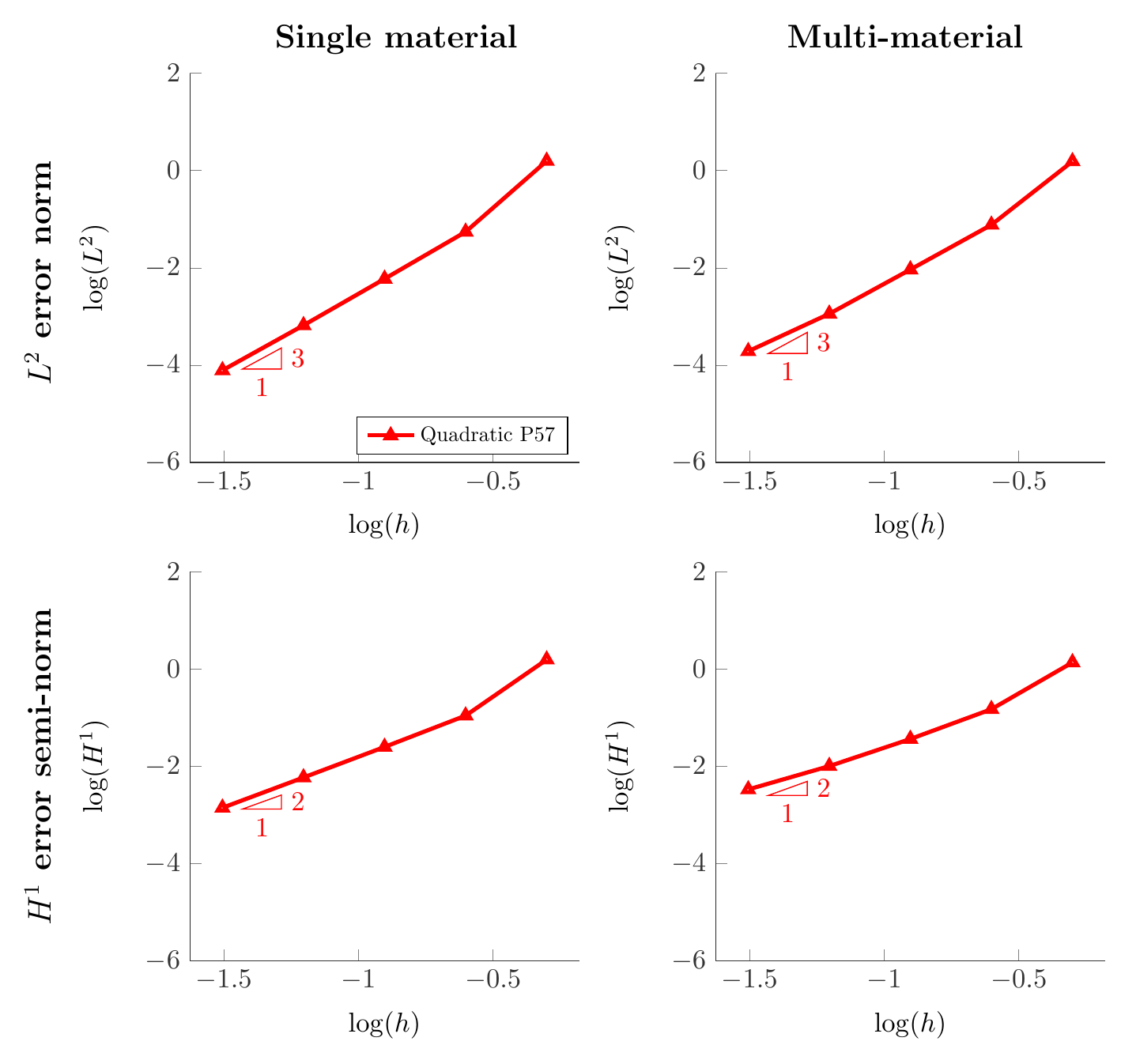}
\caption{Accuracy study on a heated single or multi-material medium in three dimensions.}
\label{Ex5_fig5}
\end{figure*}

The relative $L^2$ error norm and the relative $H^1$ error semi-norm with respect to the \textit{h}-refinement are presented in Fig. \ref{Ex5_fig5}. The highest level of refinement leads to a problem with around 15,400,000 integration elements and around 3,000,000 DOFs. Similarly to the 2D results, optimal convergence rates in the $L^2$ and $H^1$ error norms are recovered for the singe material case and close to optimal ones for the multi-material case, as singularities are present at the sharp corners of the inclusion.

This study suggests that our XIGA approach can deal with multi-material problems and provides the flexibility to enhance the approximation near singularities.

\section{Conclusions}\label{sectionConclusion}
In this paper, an XIGA approach is proposed to simulate multi-material problems. To achieve a high resolution of both the geometry and the physical response around the material interfaces, the proposed approach combines an immersed boundary technique, namely the XFEM, and the use of smooth and higher-order bases, here B-splines. Although not restricted to this approach, the geometry of the computational domains is defined using one or several LSFs in this work. The iso-levels of the LSF define subregions within the computational domain that are associated with materials. To provide additional flexibility, the computational domains are fully immersed in a background mesh. The physical responses of the systems are evaluated with the proposed XIGA approach using a novel generalized Heaviside enrichment strategy, where each basis function is enriched separately based on the material layout within the basis function support. Boundary and interface conditions are imposed weakly through Nitsche's formulation. Instabilities related to the creation of small material integration subdomains are mitigated using a generalized formulation of the face-oriented ghost penalty stabilization methodology adapted to our enrichment strategy.

The performance and the versatility of the proposed XIGA approach are studied through numerical experiments. In particular, canonical heat conduction and linear elastic problems in two and three dimensions are considered. The stability, the accuracy, and the robustness of the evaluated solutions are measured with the relative $L^2$ error norm, the relative $H^1$ error semi-norm, and the condition number.

Numerical experiments show that accurate solutions with optimal convergence rate with \textit{h}-refinement can be recovered. Additionally, the ghost penalty stabilization mitigates the effect of small material subdomains on the conditioning of the system, but also on the achieved accuracy of the solution. This remark holds regardless of the order of the B-spline basis functions, i.e., linear, quadratic, and cubic. Additionally, by studying the effect of the size of the created material integration subdomains and of the value of the ghost penalty parameter, a valid range for the ghost penalty parameter can be suggested for the considered types of problems.

Numerical examples study the application of our XIGA approach to problems with planar and curved interfaces. A slight degradation of the condition number and the solution accuracy is observed when increasing the geometric complexity. Optimal convergence rates with mesh refinement are observed for the $L^2$ and $H^1$ errors using linear, quadratic, and cubic B-splines. It should be noted that, when dealing with curved interfaces, the resolution of the geometry limits the accuracy of the evaluated state field. Individual refinement of the geometry and solution field can be used to mitigate this issue. Alternatively one can use a curved integration mesh, see \cite{2016StavrevEtAl}, or Green's theorem based integration schemes, see \cite{2022Saye}.

Moreover, numerical simulations demonstrate that \textit{N}-phase and \textit{N}-material problems can be resolved accurately. The framework allows for a large flexibility in handling and assigning phases and materials to different subregions of the computational domain. For both multi-phase and multi-material problems, optimal convergence rates with mesh refinement are recovered for the $L^2$ and $H^1$ errors using linear, quadratic, and cubic B-splines.

Several follow-up research developments are foreseen. In particular, the XIGA approach can be studied for other classes of problems, such as nonlinear and multi-physics problems. Furthermore, hierarchical B-splines, LR-splines, or T-splines can be used in place of B-splines to carry out local mesh refinement, see \cite{2012SchillingerEtAl} or \cite{2018GarauVasquez}. This procedure allows for further enhancement of the resolution of both the geometry and the physics around the interface and the accurate modelling of systems across several length scales.


\section*{Acknowledgements} 
L. No{\"e}l, M. Schmidt, J.A. Evans, and K. Maute received the support for this work from the Defense Advanced Research Projects Agency (DARPA) under the TRADES program (agreement HR0011-17-2-0022). K. Doble was supported by SANDIA through the contract PO 2120843. J. A. Evans and K. Maute were partially supported by the National Science Foundation under Grant OAC-2104106.






\begin{appendices}
\end{appendices}


\bibliographystyle{abbrvnat}
\bibliography{arXiv_LNoelEtAl2022.bib}

\end{document}